\documentclass[11pt]{article}%
\usepackage{amsmath}
\usepackage{amsfonts}
\usepackage{amssymb}
\usepackage{graphicx}
\usepackage{theorem}
\usepackage{makeidx}
\usepackage{verbatim}
\usepackage{color}%
\providecommand{\U}[1]{\protect\rule{.1in}{.1in}}
\newtheorem{thm}{Theorem}[section]

\newtheorem{pr}[thm]{Proposition}
\newtheorem{df}[thm]{Definition}
\newtheorem{rmk}[thm]{Remark}

{\theorembodyfont{\upshape}

}
\numberwithin{equation}{section} 
\begin{document}

\title{ }

\begin{center}
\vspace*{1.3cm}

\textbf{METRIC\ REGULARITY\ OF\ COMPOSITION\ SET-VALUED\ MAPPINGS: }

\textbf{METRIC\ SETTING AND\ CODERIVATIVE\ CONDITIONS}

\bigskip

by

\bigskip

MARIUS DUREA\footnote{{\small {Faculty of Mathematics, "Al. I. Cuza"
University,} \ {Bd. Carol I, nr. 11, 700506 -- Ia\c{s}i, Romania, e-mail:
\texttt{durea@uaic.ro}}}}, VAN NGAI HUYNH\footnote{{\small {Department of
Mathematics, University of Quy Nhon, 170 An Duong Vuong, Quy Nhon, Viet Nam,
e-mail: \texttt{nghiakhiem@yahoo.com}}}}, HUU TRON
NGUYEN\footnote{{\small University of Quy Nhon{, 170 An Duong Vuong, Quy Nhon,
Viet Nam,} and Laboratoire XLIM, UMR-CNRS 6172, Universit\'{e} de Limoges,
France, e-mail :{ \texttt{tron.nguyen@etu.unilim.fr}}}}, and
RADU\ STRUGARIU\footnote{{\small Department of Mathematics, "Gh. Asachi"
Technical University, \ {Bd. Carol I, nr. 11, 700506 -- Ia\c{s}i, Romania,
e-mail: \texttt{rstrugariu@tuiasi.ro}}}}

\end{center}

\bigskip

\bigskip

\noindent{\small {\textbf{Abstract:}} The paper concerns a new method to
obtain a direct proof of the openness at linear rate/metric regularity of
composite set-valued maps on metric spaces by the unification and refinement
of several methods developed somehow separately in several works of the
authors. In fact, this work is a synthesis and a precise specialization to a
general situation of some techniques explored in the last years in the
literature. In turn, these techniques are based on several important concepts
(like error bounds, lower semicontinuous envelope of a set-valued map, local
composition stability of multifunctions) and allow us to obtain two new proofs
of a recent result having deep roots in the topic of regularity of mappings.
Moreover, we make clear the idea that it is possible to use (co)derivative
conditions as tools of proof for openness results in very general situations.}

\bigskip


\bigskip


\section{Introduction}

The property of metric regularity has its origins in the open mapping
principle for linear operators obtained in the 1930s by Banach and Schauder,
and is one of the three basic and crucial principles of functional analysis,
having various applications in many branches of mathematics. Later on, this
principle was reinterpreted and generalized in two classical results: the
tangent space of Lyusternik (\cite{Lyu}) and the surjection theorem of Graves
(\cite{Gra}). The next decisive step in this history was the extension of the
Banach-Schauder principle to the case of set-valued maps with closed and
convex graph given independently by Ursescu in \cite{Urs1975} and Robinson in
\cite{Rob1976} (the celebrated Robinson-Ursescu Theorem). Moreover, it was
observed in Dmitruk, Milyutin, and Osmolovsky \cite{DMO} that the original
proof of Lyusternik from \cite{Lyu} is applicable to a much more general
setting: the sum of a covering at a rate $a>0$ single-valued mapping and a
Lipschitz one with a Lipschitz constant $b<a$ is covering at the rate $a-b$.
Another remarkable insight given in the mentioned paper is that it clearly
emphasizes the metric nature of openness and regularity properties.
Afterwards, in 1996, Ursescu \cite{Urs1996} was the first to obtain a fully
set-valued extension of the above results, in the setting of Banach spaces. On
the tracks of \cite{DMO}, the important work of Ioffe \cite{Ioffe2000} made
the crucial observation that the Lyusternik iteration process can be
successfully used when the income space is a complete metric space and the
outcome space has a linear structure with shift-invariant metric. Detailed
studies on the case of the sum of a metrically regular set-valued mapping and
a single-valued Lipschitz map appear, more recently, in works by Dontchev,
Lewis, and Rockafellar (\cite{DLR}), Dontchev, and Lewis (\cite{DL}),
Arutyunov (\cite{Arut2007}, \cite{Arut2009}), Mordukhovich (\cite{Mor2006}).
For a detailed account for the whole topic of regularity properties of
mappings, as well as various applications the reader is referred to the books
or works of many researchers: \cite{RefAz}, \cite{RefBonS}, \cite{BD},
\cite{[10]}, \cite{RefBorZu}, \cite{RefBorZ}, \cite{[17]}, \cite{RefCom},
\cite{[18]}, \cite{DLR}, \cite{RefIo3}, \cite{Ioffe2000}, \cite{Ioffe2010b},
\cite{RefJT}, \cite{JT1}, \cite{RefJAus}, \cite{Jo-Th1}, \cite{Jo-Th},
\cite{JouT1}, \cite{1JT2}, \cite{Lyu}, \cite{Mor2006}, \cite{RefMorS},
\cite{RefMorW}, \cite{NT2001}, \cite{RefNT3}, \cite{[69]}, \cite{RefPe},
\cite{Schiro}, \cite{[74]}, \cite{Zng}.

In the last years, the study of openness at linear rate (or equivalently
metric regularity) of multifunctions obtained as operations with set-valued
maps has received a new impetus coming from at least three connected issues:
the link between Lyusternik-Graves type theorems and fixed point assertions
(\cite{Arut2007}, \cite{DonFra2010}), the growing interest to generalized
forms of compositions (\cite{Ioffe2010b}, \cite{DurStr5}) and the new
developments of metric regularity results obtained under assumptions based on
generalized differentiation calculus and especially on coderivative conditions
(\cite{NTT}, \cite{DNS2011}).

Included in this stream, the present paper concerns a new method to obtain a
direct proof of the openness at linear rate/metric regularity of composite
set-valued maps on metric spaces by the unification and refinement of several
methods developed somehow separately in several works of the authors:
\cite{NT2001}, \cite{NT2008}, \cite{NTT}, \cite{NTT2}, \cite{DurStr1},
\cite{DNS2011}, \cite{DurStr5}, \cite{DurStr2012}. In some sense, this is a
synthesis and a precise specialization to a general situation of some
techniques explored in the quoted papers. In turn these techniques are based
on several important concepts (like error bounds, strong slope associated to a
function, lower semicontinuous envelope of a set-valued map, local composition
stability of multifunctions) and allow us to obtain two new proofs of a recent
result having deep roots in the literature on the topic of regularity of
mappings. Moreover, we make clear the idea that it is possible to use
(co)derivative conditions as tool of proof for openness results in very
general situations.

More precisely, the corner stones that this work rely on are mainly the
results in \cite{NT2001} on the error bounds for a nonlinear variational
system, the main result in \cite{DurStr2012} concerning the openness of a
composite multifunction, and also the coderivative conditions for metric
regularity as these appear in \cite{NTT}, \cite{DNS2011}.

The paper contains a main result (Theorem \ref{mainresult}) which is prepared
by several propositions having their own interest. In that main result one
obtains, under some already standard (hence expected) assumptions (see
\cite{DurStr5}, \cite{DurStr2012}), the openness of an auxilliary
multifunction associated to a composition set-valued map, and on this basis, a
result of openness around the reference point for the considered composition.
We want to emphasize here two main points both of them revealing the novelty
and the relevance of our work. Firstly, the conclusion is significantly richer
that the corresponding conclusions of the main results in \cite{NTT2} (from
point of view of the generality of set-valued operations), \cite{DurStr5},
\cite{DurStr2012} (from point of view of the type of openness). Secondly, the
proof is obtained using Ekeland Variational Principle (EVP, for short), a fact
that answers a question A. Ioffe raised in a discussion with the fourth author
of this paper: how to get direct proofs for openness results (and, also, for
coincidence/fixed points results), on complete metric spaces, using EVP and
not arguing by contradiction. In our knowlelge (see, for instance,
\cite{Urs1996}, \cite{DonFra2010}, \cite{DurStr5}), in many cases, the proofs
relying on EVP are made on normed vector spaces, and reasoning by
contradiction. The supplemental structure of the space (i.e., its linear
structure, but also the norm), which seems at first glance a little
surprising, it is used essentially in the construction of the contradiction.
In this work, by the careful analysis of some ideas spread in different
articles (see \cite{NT2008}, \cite{NTT}, \cite{DNS2011}), we reached at the
conclusion that the answer to the problem raised by Ioffe was already there,
but not so easy to detect: one must appropriately apply EVP for the lower
semicontinuous envelope of a certain distance function, by the appropriate
choice of an auxiliary multifunction involved in the construction of this
envelope. As consequence, by combining and extending some techniques from the
quoted articles, we are able to give here a complete and direct proof based on
EVP for the metric regularity of set-valued mappings of composition type.
Moreover, in this way, we bring more light on the links between several tools
used in getting regularity results for multifunctions.

As a by-product of the main result, a coincidence/fixed points assertion is
obtained, a fact that contributes to a discussion on this subject initiated in
\cite{Arut2007}, \cite{Arut2009}, \cite{AAGDO} and continued in
\cite{DonFra2010}, \cite{DonFra2}, \cite{DurStr5}. Futhermore, the important
role of the assertions before the main result is again emphasized, as the
(immediate) proof of the fixed point assertion relies on the appropriate
application of one of them and of the main result (Theorem \ref{mainresult}).

The last section deals with coderivative conditions for openness of composite
mapping. Here we reconsider several ideas in \cite{NTT2} and \cite{DurStr4}
from a higher point of view and we employ a calculus rule for the Fr\'{e}chet
normal cone to a intersection of sets, passing through the concept of
alliedness introduced and studied by Penot and his coauthors (\cite{Pen1998},
\cite{LPX2011}). Finally, as an interesting fact which makes the link to the
section before, we prove that one can obtain, on Asplund spaces, the
conclusion of the main result of the paper by the use of the coderivative
condition previously developed.

\section{Preliminaries}

This section contains some basic definitions and results used in the sequel.
In what follows, we suppose that the involved spaces are metric spaces, unless
otherwise stated. In this setting, $B(x,r)$ and $\overline{B}(x,r)$ denote the
open and the closed ball with center $x$ and radius $r,$ respectively. On a
product space we usually take the additive metric; when we choose another
metric, this will be stated explicitly. If $x\in X$ and $A\subset X,$ one
defines the distance from $x$ to $A$ as $d(x,A):=\inf\{d(x,a)\mid a\in A\}.$
As usual, we use the convention $d(x,\emptyset)=\infty.$ The excess from a set
$A$ to a set $B$ is defined as $e(A,B):=\sup\{d(a,B)\mid a\in A\}.$ For a
non-empty set $A\subset X$ we put $\operatorname*{cl}A$ for its topological
closure. One says that a set $A$ is locally complete (closed) if there exists
$r>0$ such that $A\cap\overline{B}(x,r)$ is complete (closed). The symbol
$\mathcal{V}(x)$ stands for the system of the neighborhoods of $x.$

Let $X,Y,Z,P$ be metric spaces. For a multifunction $F:X\rightrightarrows Y,$
the graph of $F$ is the set $\operatorname*{Gr}F:=\{(x,y)\in X\times Y\mid
y\in F(x)\}.$ If $A\subset X,$ then $F(A):=%
{\displaystyle\bigcup\limits_{x\in A}}
F(x).$ The inverse set-valued map of $F$ is $F^{-1}:Y\rightrightarrows X$
given by $F^{-1}(y)=\{x\in X\mid y\in F(x)\}$. If $F_{1}:X\rightrightarrows
Y,F_{2}:X\rightrightarrows Z,$ we define the set-valued map $(F_{1}%
,F_{2}):X\rightrightarrows Y\times Z$ by $(F_{1},F_{2})(x):=F_{1}(x)\times
F_{2}(x).$ For a parametric multifunction $F:X\times P\rightrightarrows Y,$ we
use the notations: $F_{p}(\cdot):=F(\cdot,p)$ and $F_{x}(\cdot):=F(x,\cdot).$

\bigskip

We recall now the concepts of openness at linear rate, metric regularity and
Aubin property of a multifunction around the reference point.

\begin{df}
\label{around}Let $F:X\rightrightarrows Y$ be a multifunction and
$(\overline{x},\overline{y})\in\operatorname{Gr}F.$

(i) $F$ is said to be open at linear rate $L>0$ (or $L-$open) around
$(\overline{x},\overline{y})$ if there exist a positive number $\varepsilon>0$
and two neighborhoods $U\in\mathcal{V}(\overline{x}),$ $V\in\mathcal{V}%
(\overline{y})$ such that, for every $\rho\in(0,\varepsilon)$ and every
$(x,y)\in\operatorname*{Gr}F\cap\lbrack U\times V],$%
\begin{equation}
B(y,\rho L)\subset F(B(x,\rho)). \label{Lopen}%
\end{equation}

The supremum of $L>0$ over all the combinations $(L,U,V,\varepsilon)$ for
which (\ref{Lopen}) holds is denoted by $\operatorname*{lop}F(\overline
{x},\overline{y})$ and is called the exact linear openness bound, or the exact
covering bound of $F$ around $(\overline{x},\overline{y}).$

(ii) $F$ is said to have the Aubin property (or to be Lipschitz-like) around
$(\overline{x},\overline{y})$ with constant $L>0$ if there exist two
neighborhoods $U\in\mathcal{V}(\overline{x}),$ $V\in\mathcal{V}(\overline{y})$
such that, for every $x,u\in U,$%
\begin{equation}
e(F(x)\cap V,F(u))\leq Ld(x,u). \label{LLip_like}%
\end{equation}

The infimum of $L>0$ over all the combinations $(L,U,V)$ for which
(\ref{LLip_like}) holds is denoted by $\operatorname*{lip}F(\overline
{x},\overline{y})$ and is called the exact Lipschitz bound of $F$ around
$(\overline{x},\overline{y}).$

(iii) $F$ is said to be metrically regular with constant $L>0$ around
$(\overline{x},\overline{y})$ if there exist two neighborhoods $U\in
\mathcal{V}(\overline{x}),$ $V\in\mathcal{V}(\overline{y})$ such that, for
every $(x,y)\in U\times V,$%
\begin{equation}
d(x,F^{-1}(y))\leq Ld(y,F(x)). \label{Lmreg}%
\end{equation}

The infimum of $L>0$ over all the combinations $(L,U,V)$ for which
(\ref{Lmreg}) holds is denoted by $\operatorname*{reg}F(\overline{x}%
,\overline{y})$ and is called the exact regularity bound of $F$ around
$(\overline{x},\overline{y}).$
\end{df}

The links between the previous notions are as follows (see, e.g.,
\cite[Theorem 9.43]{RocWet}, \cite[Theorems 1.52]{Mor2006}, \cite[Section
3E]{DontRock2009b}).

\begin{thm}
\label{link_around}Let $F:X\rightrightarrows Y$ be a multifunction and
$(\overline{x},\overline{y})\in\operatorname{Gr}F.$ Then $F$ is open at linear
rate around $(\overline{x},\overline{y})$ iff $F^{-1}$ has the Aubin property
around $(\overline{y},\overline{x})$ iff $F$ is metrically regular around
$(\overline{x},\overline{y})$. Moreover, in every of the previous situations,%
\[
(\operatorname*{lop}F(\overline{x},\overline{y}))^{-1}=\operatorname*{lip}%
F^{-1}(\overline{y},\overline{x})=\operatorname*{reg}F(\overline{x}%
,\overline{y}).
\]

\end{thm}

In the case of parametric set-valued maps one has the following partial
notions of linear openness, metric regularity and Aubin property around the
reference point.

\begin{df}
Let $F:X\times P\rightrightarrows Y$ be a multifunction and $((\overline
{x},\overline{p}),\overline{y})\in\operatorname{Gr}F.$

(i) $F$ is said to be open at linear rate $L>0$ (or $L-$open) with respect to
$x$ uniformly in $p$ around $((\overline{x},\overline{p}),\overline{y})$ if
there exist a positive number $\varepsilon>0$ and some neighborhoods
$U\in\mathcal{V}(\overline{x}),$ $V\in\mathcal{V}(\overline{p}),$
$W\in\mathcal{V}(\overline{y})$ such that, for every $\rho\in(0,\varepsilon),$
every $p\in V$ and every $(x,y)\in\operatorname*{Gr}F_{p}\cap\lbrack U\times
W],$%
\begin{equation}
B(y,\rho L)\subset F_{p}(B(x,\rho)). \label{pLopen}%
\end{equation}

The supremum of $L>0$ over all the combinations $(L,U,V,W,\varepsilon)$ for
which (\ref{pLopen}) holds is denoted by $\widehat{\operatorname*{lop}}%
_{x}F((\overline{x},\overline{p}),\overline{y})$ and is called the exact
linear openness bound, or the exact covering bound of $F$ in $x$ around
$((\overline{x},\overline{p}),\overline{y}).$

(ii) $F$ is said to have the Aubin property (or to be Lipschitz-like) with
respect to $x$ uniformly in $p$ around $((\overline{x},\overline{p}%
),\overline{y})$ with constant $L>0$ if there exist some neighborhoods
$U\in\mathcal{V}(\overline{x}),$ $V\in\mathcal{V}(\overline{p}),$
$W\in\mathcal{V}(\overline{y})$ such that, for every $x,u\in U$ and every
$p\in V,$%
\begin{equation}
e(F_{p}(x)\cap W,F_{p}(u))\leq Ld(x,u). \label{pLLip_like}%
\end{equation}

The infimum of $L>0$ over all the combinations $(L,U,V,W)$ for which
(\ref{pLLip_like}) holds is denoted by $\widehat{\operatorname*{lip}}%
_{x}F((\overline{x},\overline{p}),\overline{y})$ and is called the exact
Lipschitz bound of $F$ in $x$ around $((\overline{x},\overline{p}%
),\overline{y}).$

(iii) $F$ is said to be metrically regular with constant $L>0$ with respect to
$x$ uniformly in $p$ around $((\overline{x},\overline{p}),\overline{y})$ if
there exist some neighborhoods $U\in\mathcal{V}(\overline{x}),$ $V\in
\mathcal{V}(\overline{p}),$ $W\in\mathcal{V}(\overline{y})$ such that, for
every $(x,p,y)\in U\times V\times W,$%
\begin{equation}
d(x,F_{p}^{-1}(y))\leq Ld(y,F_{p}(x)). \label{pLmreg}%
\end{equation}

The infimum of $L>0$ over all the combinations $(L,U,V,W)$ for which
(\ref{pLmreg}) holds is denoted by $\widehat{\operatorname*{reg}}%
_{x}F((\overline{x},\overline{p}),\overline{y})$ and is called the exact
regularity bound of $F$ in $x$ around $((\overline{x},\overline{p}%
),\overline{y}).$
\end{df}

Interchanging the roles of $p$ and $x$ one gets a similar set of concepts.

\bigskip

Let now $X$ be a metric space, and $f:X\rightarrow\mathbb{R\cup\{+}\infty\}$
be a function. Set%

\[
S:=\{x\in X:f(x)\leq0\}.
\]

One denotes the quantity $\max\{f(x),0\}$ by $[f(x)]_{+}.$

Here and in what follows the convention $0\cdot(+\infty)=0$ is used. The next
result gives an error bound estimation and it will be useful as an
intermediate step in the proof of the main result of the paper.

\begin{thm}
\label{m}(\cite[Theorem 2.1]{NT2008}) Let $(X,d)$ be a complete metric space,
and let $f:X\rightarrow\mathbb{R\cup\{+}\infty\}$ be a lower semicontinuous
function, and $\overline{x}\in X$ such that $f(\overline{x})>0.$ Setting%
\[
m(\overline{x}):=\inf\left\{  \sup_{u\in X,u\not =\overline{x}}\frac
{f(x)-[f(u)]_{+}}{d(x,u)}:%
\begin{array}
[c]{c}%
d(x,\overline{x})<d(\overline{x},S)\\
f(x)\leq f(\overline{x})
\end{array}
\right\}  ,
\]

\noindent one has%
\[
m(\overline{x})\cdot d(\overline{x},S)\leq f(\overline{x}).
\]

\end{thm}

In what follows we shall use the lower semicontinuous envelope associated to a
multifunction $F:X\rightrightarrows Y$ as $\varphi_{F}:X\times Y\rightarrow
\mathbb{R\cup\{+}\infty\},$%
\[
\varphi_{F}(x,y):=\liminf_{(u,v)\rightarrow(x,y)}d(v,F(u))=\liminf
_{u\rightarrow x}d(y,F(u)).
\]

This concept, with appropriate choices of the base multifunction, proved to be
an useful tool in proving results on metric regularity of set-valued mappings
and in vector optimization, using the error bound approach (see, e.g.,
\cite{NTT}, \cite{DNS2011}). Observe that $\varphi_{F}(x,y)\geq0$ and
$\varphi_{F}(x,y)\leq d(y,F(x))$ for every $(x,y)\in X\times Y.$

\begin{pr}
\label{F-1}(\cite[Lemma 1]{DNS2011}) Let $F:X\rightrightarrows Y$ be a
multifunction. Then, for every $y\in Y,$%
\[
\{x\in X:\varphi_{F}(x,y)=0\}=(\operatorname*{cl}F)^{-1}(y),
\]
where $\operatorname*{cl}F$ is the multifunction whose graph is
$\operatorname*{cl}\operatorname*{Gr}F.$ In particular, if $F$ has closed
graph, then%
\[
F^{-1}(y)=\{x\in X:\varphi_{F}(x,y)=0\}.
\]

\end{pr}

\section{Metric regularity of compositions}

Let $F_{1}:X\rightrightarrows Y_{1}$, $F_{2}:X\rightrightarrows Y_{2}$ and
$G:Y_{1}\times Y_{2}\rightrightarrows Z$ be set-valued mappings, where
$X,Y_{1},Y_{2},Z$ are metric spaces. Consider the following composition
multifunctions $H:X\rightrightarrows Z$ defined by
\begin{equation}
H(x):=G(F_{1}(x),F_{2}(x)). \label{compomap}%
\end{equation}
Our aim in this article is to investigate the metric regularity of this
multifunction $H$. As in the case of the sum of multifunctions (that is,
$Y_{1}=Y_{2}=Y,$ where $Y$ is a normed vector space, and $G(y_{1},y_{2}%
)=y_{1}+y_{2}$), in general, $H$ is not necessarily closed. Thus almost all
the results on metric regularity known in the literature, that need the
closedness of the multifunction under consideration, could not be directly
applied for $H$. For this purpose, let us consider, associated to $F_{1}%
,F_{2}$ and $G$ the following multifunction: $R:X\times Y_{1}\times
Y_{2}\rightrightarrows Z,$ given as
\begin{equation}
R(x,y_{1},y_{2})=\left\{
\begin{array}
[c]{ll}%
G(y_{1},y_{2}), & \text{if }(y_{1},y_{2})\in F_{1}(x)\times F_{2}(x)\\
\emptyset, & \text{otherwise.}%
\end{array}
\right.  \label{R}%
\end{equation}

Obviously, $R$ is a closed-graph multifunction whenever $F_{1},F_{2},G$ are
closed-graph ones.

Next, we want to see how the associated function $\varphi_{R}$ looks like. We
know by definition that
\[
\varphi_{R}((x,y_{1},y_{2}),z):=\liminf_{(x^{\prime},y_{1}^{\prime}%
,y_{2}^{\prime},z^{\prime})\rightarrow(x,y_{1},y_{2},z)}d(z^{\prime
},R(x^{\prime},y_{1}^{\prime},y_{2}^{\prime})).
\]
Clearly, if $(y_{1},y_{2})\notin F_{1}(x)\times F_{2}(x),$ say $y_{1}%
\not \in F_{1}(x),$ and $F_{1}(x)$ is closed, then $\varphi_{R}((x,y_{1}%
,y_{2}),z)=+\infty.$ As consequence, if $F_{1},F_{2}$ are closed-graph
multifunctions and $(y_{1},y_{2})\notin F_{1}(x)\times F_{2}(x),$ one has that
$\varphi_{R}((x,y_{1},y_{2}),z)=+\infty.$

Otherwise, we have
\begin{align*}
\varphi_{R}((x,y_{1},y_{2}),z)  &  ={{\liminf\limits_{\substack{(x^{\prime
},y_{1}^{\prime},y_{2}^{\prime},z^{\prime})\rightarrow(x,y_{1},y_{2}%
,z)\\{(y_{1}^{\prime},y_{2}^{\prime})\in F_{1}(x^{\prime})\times
F_{2}(x^{\prime})}}}}}d(z^{\prime},R(x^{\prime},y_{1}^{\prime},y_{2}^{\prime
}))\\
&  ={{\liminf\limits_{\substack{{(x^{\prime},y_{1}^{\prime},y_{2}^{\prime
})\rightarrow(x,y_{1},y_{2})}\\{(y_{1}^{\prime},y_{2}^{\prime})\in
F_{1}(x^{\prime})\times F_{2}(x^{\prime})}}}}}d(z,R(x^{\prime},y_{1}^{\prime
},y_{2}^{\prime}))\\
&  ={{\liminf\limits_{\substack{{(x^{\prime},y_{1}^{\prime},y_{2}^{\prime
})\rightarrow(x,y_{1},y_{2})}\\{(y_{1}^{\prime},y_{2}^{\prime})\in
F_{1}(x^{\prime})\times F_{2}(x^{\prime})}}}}}d(z,G(y_{1}^{\prime}%
,y_{2}^{\prime})).
\end{align*}

In conclusion, if $F_{1},F_{2}$ are closed-graph multifunctions, one has that%
\[
\varphi_{R}((x,y_{1},y_{2}),z)=\left\{
\begin{array}
[c]{ll}%
{{\liminf\limits_{\substack{{(x^{\prime},y_{1}^{\prime},y_{2}^{\prime
})\rightarrow(x,y_{1},y_{2})}\\{(y_{1}^{\prime},y_{2}^{\prime})\in
F_{1}(x^{\prime})\times F_{2}(x^{\prime})}}}}}d(z,G(y_{1}^{\prime}%
,y_{2}^{\prime})), & \text{if }(y_{1},y_{2})\in F_{1}(x)\times F_{2}(x)\\
+\infty, & \text{otherwise.}%
\end{array}
\right.
\]

We are ready to prove a result which will become very useful in the sequel. As
mentioned, in the next result we take the additive distance on product spaces.

\begin{pr}
\label{proMT} Let $F_{1}:X\rightrightarrows Y_{1}$, $F_{2}:X\rightrightarrows
Y_{2}$ and $G:Y_{1}\times Y_{2}\rightrightarrows Z$ be multifunctions and
$\overline{z}\in G(\overline{y}_{1},\overline{y}_{2}),(\overline{y}%
_{1},\overline{y}_{2})\in F_{1}(\overline{x})\times F_{2}(\overline{x})$.
Consider the following statements:

{(i)} there exist a neighborhood $\mathcal{U}\times\mathcal{V}_{1}%
\times\mathcal{V}_{2}\times\mathcal{W}\subset X\times Y_{1}\times Y_{2}\times
Z$ of $(\overline{x},\overline{y}_{1},\overline{y}_{2},\overline{z})$ and
$\tau>0$ such that
\begin{equation}
d((x,y_{1},y_{2}),R^{-1}(z))\leq\tau\varphi_{R}((x,y_{1},y_{2}),z)\quad
\text{for all }(x,y_{1},y_{2},z)\in\mathcal{U}\times\mathcal{V}_{1}%
\times\mathcal{V}_{2}\times\mathcal{W}. \label{i}%
\end{equation}

{(ii)} there exist a neighborhood \ $\mathcal{U}\times\mathcal{V}_{1}%
\times\mathcal{V}_{2}\times\mathcal{W}\subset X\times Y_{1}\times Y_{2}\times
Z$ of $(\overline{x},\overline{y}_{1},\overline{y}_{2},\overline{z})$ and
$\tau>0$ such that
\begin{equation}
d(x,H^{-1}(z))\leq\tau d(z,G(F_{1}(x)\cap\mathcal{V}_{1},F_{2}(x)\cap
\mathcal{V}_{2}))\quad\text{for all }(x,z)\in\mathcal{U}\times\mathcal{W}.
\label{ii}%
\end{equation}

{(iii)} there exists $\varepsilon>0$ such that for every $\rho\in
(0,\varepsilon)$ and for every $(x,y_{1},y_{2},z)\in B(\overline
{x},\varepsilon)\times B(\overline{y}_{1},\varepsilon)\times B(\overline
{y}_{2},\varepsilon)\times B(\overline{z},\varepsilon)$ such that $z\in
G(y_{1},y_{2}),(y_{1},y_{2})\in F_{1}(x)\times F_{2}(x),$
\[
B(z,\rho\tau^{-1})\in H(B(x,\rho)).
\]
Then, we have the following implications: $(i)\Rightarrow(ii)\Rightarrow
(iii).$
\end{pr}

\noindent\textbf{Proof.} $(i)\Rightarrow(ii).$ Using $(i)$, we get that there
exist a neighborhood $\mathcal{U}\times\mathcal{V}_{1}\times\mathcal{V}%
_{2}\times\mathcal{W}\subset X\times Y_{1}\times Y_{2}\times Z$ of
$(\overline{x},\overline{y}_{1},\overline{y}_{2},\overline{z})$ and $\tau>0$
such that%
\begin{equation}
d((x,y_{1},y_{2}),R^{-1}(z))\leq\tau d(z,G(y_{1},y_{2})) \label{regR}%
\end{equation}
for all $(x,y_{1},y_{2},z)\in\mathcal{U}\times\mathcal{V}_{1}\times
\mathcal{V}_{2}\times\mathcal{W}\;$with$\;(y_{1},y_{2})\in F_{1}(x)\times
F_{2}(x).$

Take $z\in\mathcal{W}$. If for $x\in\mathcal{U},$ $F_{1}(x)\cap\mathcal{V}%
_{1}=\emptyset$ or $F_{2}(x)\cap\mathcal{V}_{2}=\emptyset,$ then
$d(z,G(F_{1}(x)\cap\mathcal{V}_{1},F_{2}(x)\cap\mathcal{V}_{2}))=+\infty.$
Suppose now that there is $(y_{1},y_{2})\in\lbrack F_{1}(x)\cap\mathcal{V}%
_{1}]\times\lbrack F_{2}(x)\cap\mathcal{V}_{2}]$ such that $G(y_{1}%
,y_{2})\not =\emptyset,$ because otherwise, again, relation (\ref{ii})
trivially holds.

Take now arbitrary $(y_{1},y_{2})\in\lbrack F_{1}(x)\cap\mathcal{V}_{1}%
]\times\lbrack F_{2}(x)\cap\mathcal{V}_{2}]$ such that $G(y_{1},y_{2}%
)\not =\emptyset.$ From (\ref{regR}), it follows that for every $\varepsilon
>0,$ there exists $(x^{\prime},y_{1}^{\prime},y_{2}^{\prime})\in R^{-1}(z)$,
i.e., $z\in G(y_{1}^{\prime},y_{2}^{\prime}),$ $(y_{1}^{\prime},y_{2}^{\prime
})\in F_{1}(x^{\prime})\times F_{2}(x^{\prime}),$ such that
\begin{equation}
d((x,y_{1},y_{2}),(x^{\prime},y_{1}^{\prime},y_{2}^{\prime}))\leq\tau
d(z,G(y_{1},y_{2}))+\varepsilon. \label{ineg_dist1}%
\end{equation}
Consequently,
\begin{equation}
d(x,x^{\prime})\leq\tau d(z,G(y_{1},y_{2}))+\varepsilon, \label{ineg_dist2}%
\end{equation}
and, since $z\in G(y_{1}^{\prime},y_{2}^{\prime}),$ $(y_{1}^{\prime}%
,y_{2}^{\prime})\in F_{1}(x^{\prime})\times F_{2}(x^{\prime})$, we have that
$x^{\prime}\in H^{-1}(z),$ hence%
\[
d(x,H^{-1}(z))\leq\tau d(z,G(y_{1},y_{2}))+\varepsilon,
\]

\noindent and, finally,%
\[
d(x,H^{-1}(z))\leq\tau d(z,G(F_{1}(x)\cap\mathcal{V}_{1},F_{2}(x)\cap
\mathcal{V}_{2}))+\varepsilon\quad\text{for all }(x,z)\in\mathcal{U}%
\times\mathcal{W}.
\]

Making $\varepsilon\rightarrow0,$ one gets the conclusion.

$(ii)\Rightarrow(iii).$ By $(ii)$, there are $\delta,\tau>0$ such that for
every $(x,y_{1},y_{2},z)\in B(\overline{x},\delta)\times B(\overline{y}%
_{1},\delta)\times B(\overline{y}_{2},\delta)\times B(\overline{z},\delta),$
one has
\[
d(x,H^{-1}(z))\leq\tau d(z,G(F_{1}(x)\cap B(\overline{y}_{1},\delta
),F_{2}(x)\cap B(\overline{y}_{2},\delta)))\quad\text{for all }(x,z)\in
B(\overline{x},\delta)\times B(\overline{z},\delta).
\]
Choose $\varepsilon<\delta\dfrac{\tau}{\tau+1}$. Then, for every
$(x,y_{1},y_{2},z)\in B(\overline{x},\varepsilon)\times B(\overline{y}%
_{1},\varepsilon)\times B(\overline{y}_{2},\varepsilon)\times B(\overline
{z},\varepsilon)$ such that $z\in G(y_{1},y_{2}),$ $(y_{1},y_{2})\in
F_{1}(x)\times F_{2}(x)$ and every $\rho\in(0,\varepsilon),$ take
\[
y\in B(z,\rho\tau^{-1}).
\]
Then
\[
d(y,\overline{z})\leq d(y,z)+d(z,\overline{z})<\rho\tau^{-1}+\varepsilon
<\varepsilon\tau^{-1}+\varepsilon<\frac{\tau+1}{\tau}\delta\frac{\tau}{\tau
+1}=\delta.
\]
It follows that
\[
d(x,H^{-1}(y))\leq\tau d(y,z)<\tau\rho\tau^{-1}=\rho.
\]
Thus, there is $u\in H^{-1}(y)$ or, equivalently, $y\in H(u)$ such that
$d(x,u)<\rho.$ It means that
\[
y\in H(B(x,\rho).
\]
In conclusion,
\[
B(z,\rho\tau^{-1})\in H(B(x,\rho)).
\]
The proof is complete.$\hfill\square$

\bigskip

The next concept of local composition stability, initially introduced in
\cite{DurStr2012} in order to conserve the Aubin property for compositions, as
a natural extension of its corresponding predecesor, the local sum-stability
(\cite{DurStr4}), will allow us to obtain regularity results around the
reference point.

\begin{df}
Let $F:X\rightrightarrows Y,$ $G:Y\rightrightarrows Z$ be multifunctions and
$(\overline{x},\overline{y},\overline{z})\in X\times Y\times Z$ such that
$\overline{y}\in F(\overline{x}),$ $\overline{z}\in G(\overline{y}).$ We say
that the multifunctions $F,G$ are locally composition-stable around
$(\overline{x},\overline{y},\overline{z})$ if for every $\varepsilon>0$ there
exists $\delta>0$ such that, for every $x\in B(\overline{x},\delta)$ and every
$z\in(G\circ F)(x)\cap B(\overline{z},\delta),$ there is $y\in F(x)\cap
B(\overline{y},\varepsilon)$ such that $z\in G(y).$
\end{df}

A simple case which ensures the locally composition-stable property of $H$ is
as follows.\newline

\begin{pr}
\label{usc-composition} Let $X, Y_{1}, Y_{2}, Z$ be metric spaces,
$F_{1}:X\rightrightarrows Y_{1},$ $F_{2}:X\rightrightarrows Y_{2}$ and
$G:Y_{1}\times Y_{2}\rightrightarrows Z$ be three multifunctions, and
$(\overline{x},\overline{y}_{1},\overline{y}_{2},\overline{z})\in X\times
Y_{1}\times Y_{2}\times Z$ be such that $F_{1}(\overline{x})=\{\overline
{y}_{1}\},$ $F_{2}(\overline{x})=\{\overline{y}_{2}\}$ and $F_{1},F_{2}$ are
upper semicontinuous at $\overline{x}.$ Then $(F_{1},F_{2}),G$ are locally
composition-stable around $(\overline{x},(\overline{y}_{1},\overline{y}%
_{2}),\overline{z}).$
\end{pr}

\noindent\textbf{Proof.} Since $F_{1},F_{2}$ are upper semicontinuous at
$\overline{x},$ for every $\varepsilon>0$ there exist $\delta_{1},\delta
_{2}>0$ such that
\begin{align}
F_{1}(x)  &  \subset B(\overline{y}_{1},2^{-1}\varepsilon),\quad\text{for all
}x\in B(\overline{x},\delta_{1}),\label{usc}\\
F_{2}(x)  &  \subset B(\overline{y}_{2},2^{-1}\varepsilon),\quad\text{for all
}x\in B(\overline{x},\delta_{2}).\nonumber
\end{align}
Set $\delta:=\min\{\delta_{1},\delta_{2},\varepsilon\},$ and take $x\in
B(\overline{x},\delta),$ $z\in(G\circ(F_{1},F_{2}))(x)\cap B(\overline
{z},\delta)$. Then, there are $y_{1}\in F_{1}(x),$ $y_{2}\in F_{2}(x)$ such
that $z\in G(y_{1},y_{2}).$ By (\ref{usc}), one gets that $y_{1}\in
B(\overline{y}_{1},2^{-1}\varepsilon),$ $y_{2}\in B(\overline{y}_{2}%
,2^{-1}\varepsilon),$ i.e. the conclusion.$\hfill\square$

\bigskip

Another interesting case which ensures the locally composition-stable property
of $H$ is given in the following proposition.

\begin{pr}
\label{usc-compositionnew} Let $X,Y_{1},Y_{2},Z$ be metric spaces,
$F_{1}:X\rightrightarrows Y_{1},$ $F_{2}:X\rightrightarrows Y_{2}$ be two
set-valued mappings and $G:Y_{1}\times Y_{2}\rightarrow Z$ be a single-valued
mapping, and $(\overline{x},\overline{y}_{1},\overline{y}_{2},\overline{z})\in
X\times Y_{1}\times Y_{2}\times Z$ be such that $\overline{z}=G(\overline
{y}_{1},\overline{y}_{2}),(\overline{y}_{1},\overline{y}_{2})\in
F_{1}(\overline{x})\times F_{2}(\overline{x})$ and

{(i)} $F_{1}(\overline{x})=\{\overline{y}_{1}\},$ and $F_{1}$ is upper
semicontinuous at $\overline{x};$

({ii)} $G(\cdot,\overline{y}_{2})$ is continuous at $\overline{y}_{1};$

{(iii)} $G(y_{1},\cdot)$ is an isometry for all $y_{1}$ near $\overline{y}%
_{1}.$

\noindent Then $(F_{1},F_{2}),G$ are locally composition-stable around
$(\overline{x},(\overline{y}_{1},\overline{y}_{2}),\overline{z}).$
\end{pr}

\noindent\textbf{Proof.} Let $\varepsilon>0$ be given. By $(ii)$, there exists
$\delta_{1}>0$ such that
\begin{equation}
\label{Compo1}d(G(y_{1},\overline{y}_{2}),G(\overline{y}_{1},\overline{y}%
_{2}))<\varepsilon/2, \;\text{for all}\; y_{1}\in B(\overline{y}_{1}%
,\delta_{1}).
\end{equation}

By $(i)$, there is $\delta_{2}>0$ such that
\begin{equation}
F_{1}(x)\subset B(\overline{y}_{1},\delta_{1})\;\text{for all}\;x\in
B(\overline{x},\delta_{2}). \label{Compo2}%
\end{equation}
Let $\delta:=\{\delta_{1},\delta_{2},\varepsilon/2\},$ and let $x\in
B(\overline{x},\delta),$ $z\in G(F_{1}(x),F_{2}(x))\cap B(\overline{z}%
,\delta).$

Then, there are $y_{1}\in F_{1}(x),y_{2}\in F_{2}(x)$ such that $z=G(y_{1}%
,y_{2})\in B(\overline{z},\delta).$ Moreover, by (\ref{Compo2}), one has
$y_{1}\in B(\overline{y}_{1},\delta_{1})$. Thus, by (\ref{Compo1}),
$d(G(y_{1},\overline{y}_{2}),G(\overline{y}_{1},\overline{y}_{2}%
))<\varepsilon/2.$ Consequently, by $(iii)$, we obtain that
\begin{align*}
\delta &  >d(z,\overline{z})=d(G(y_{1},y_{2}),G(\overline{y}_{1},\overline
{y}_{2}))\geq\\
&  \geq d(G(y_{1},y_{2}),G({y}_{1},\overline{y}_{2}))-d(G(y_{1},\overline
{y}_{2}),G(\overline{y}_{1},\overline{y}_{2}))\geq d(y_{2},\overline{y}%
_{2})-\varepsilon/2.
\end{align*}
So,
\[
d(y_{2},\overline{y}_{2})<\delta+\varepsilon/2<\varepsilon.
\]
The proof is complete.$\hfill\square$

\bigskip

\begin{pr}
\label{MTcomsta} Let $X,Y_{1},Y_{2},Z$ be metric spaces, $F_{1}%
:X\rightrightarrows Y_{1},$ $F_{2}:X\rightrightarrows Y_{2}$ and
$G:Y_{1}\times Y_{2}\rightrightarrows Z$ be multifunctions and $(\overline
{x},\overline{y}_{1},\overline{y}_{2},\overline{z})\in X\times Y_{1}\times
Y_{2}\times Z$ be such that $\overline{z}\in G(\overline{y}_{1},\overline
{y}_{2}),(\overline{y}_{1},\overline{y}_{2})\in F_{1}(\overline{x})\times
F_{2}(\overline{x})$. If there exist $\delta>0$ and $\tau>0$ such that
\begin{equation}
d(x,H^{-1}(z))\leq\tau d(z,G(F_{1}(x)\cap B(\overline{y}_{1},\delta
),F_{2}(x)\cap B(\overline{y}_{2},\delta)))\quad\text{for all }(x,z)\in
B(\overline{x},\delta)\times B(\overline{z},\delta) \label{equ1}%
\end{equation}
and $(F_{1},F_{2}),G$ are locally composition-stable around $(\overline
{x},(\overline{y}_{1},\overline{y}_{2}),\overline{z}),$ then $H$ is metrically
regular around $(\overline{x},\overline{z})$ with modulus $\tau.$

Moreover, if $F_{1}(\overline{x})=\{\overline{y}_{1}\},$ $F_{2}(\overline
{x})=\{\overline{y}_{2}\},$ $F_{1},F_{2}$ are upper semicontinuous at
$\overline{x}$ and (\ref{equ1}) holds, then $H$ is metrically regular around
$(\bar{x},\bar{z})$ with modulus $\tau.$
\end{pr}

\noindent\textbf{Proof.} Suppose that (\ref{equ1}) holds. Since $(F_{1}%
,F_{2}),G$ are locally composition-stable around the point $(\overline
{x},(\overline{y}_{1},\overline{y}_{2}),\overline{z})$, there exists $\eta>0$
such that for every $x\in B(\overline{x},\eta)$ and every $z\in(G\circ
(F_{1},F_{2}))(x)\cap B(\overline{z},\eta),$ there is $(y_{1},y_{2})\in
(F_{1}(x)\cap B(\overline{y}_{1},\delta))\times(F_{2}(x)\cap B(\overline
{y}_{2},\delta))$ such that $z\in G(y).$

Taking $\eta$ smaller if necessary, we can assume that $\eta<\delta$. Consider
the following two cases:

Case 1. $d(z,H(x))<\eta/2.$

\noindent Fix $(x,z)\in B(\overline{x},\eta/2)\times B(\overline{z},\eta
/2)\ $and $\gamma>0$, small enough in order to have
\[
d(z,H(x))+\gamma<\eta/2.
\]
Then, there exists $t\in G(y_{1},y_{2}),(y_{1},y_{2})\in F_{1}(x)\times
F_{2}(x)$ such that $d(z,t)<d(z,G(F_{1}(x),F_{2}(x))+\gamma<\eta/2.$ Since
$d(z,\overline{z})<\eta/2$, one has that
\[
d(t,\overline{z})\leq d(z,t)+d(z,\overline{z})<\eta/2+\eta/2=\eta.
\]
Thus,
\[
t\in(G\circ(F_{1},F_{2}))(x)\cap B(\overline{z},\eta).
\]
Then, by the locally composition-stable property of $(F_{1},F_{2}),G$ around
$(\overline{x},(\overline{y}_{1},\overline{y}_{2}),\overline{z})$, there is
$(y_{1},y_{2})\in(F_{1}(x)\cap B(\overline{y}_{1},\delta))\times(F_{2}(x)\cap
B(\overline{y}_{2},\delta))$ such that $t\in G(y_{1},y_{2}).$

\noindent Consequently,
\[
d(x,H^{-1}(z))\leq\tau d(z,G(F_{1}(x)\cap B(\overline{y}_{1},\delta
),F_{2}(x)\cap B(\overline{y}_{2},\delta)))\leq\tau d(z,t)<\tau(d(z,G(F_{1}%
(x),F_{2}(x)))+\gamma).
\]
As $\gamma>0$ is arbitrarily small, one obtains that
\[
d(x,H^{-1}(z))\leq\tau d(z,G(F_{1}(x),F_{2}(x)))=\tau d(z,H(x)).
\]
Since $(x,z)$ is arbitrary in $B(\overline{x},\eta/2)\times B(\overline
{z},\eta/2),$ this yields the conclusion.

Case 2. $d(z,H(x))\geq\eta/2.$

\noindent Without loosing the generality, choose $\eta$ sufficiently small so
that $\tau\eta/2<\delta.$ Fix now $(x,z)\in B(\overline{x},\tau\eta/4)\times
B(\overline{z},\eta/4).$

If $d(\overline{x},H^{-1}(z))=0,$ one can succesively write%
\[
d(x,H^{-1}(x))\leq d(x,\overline{x})+d(\overline{x},H^{-1}(z))\leq\frac{\eta
}{2}\cdot\frac{\tau}{2}\leq\frac{\tau}{2}d(z,H(x))<\tau d(z,H(x)).
\]

\noindent Suppose now $d(\overline{x},H^{-1}(z))>0.$ Then, for any
$\varepsilon>0$, by (\ref{equ1}), there exists $u\in H^{-1}(z)$ such that
\[
d(\overline{x},u)<(1+\varepsilon)\tau d(z,H(\overline{x}))\leq(1+\varepsilon
)\tau d(z,\overline{z})<(1+\varepsilon)\tau\frac{\eta}{4}\leq(1+\varepsilon
)\frac{\tau}{2}d(z,H(x)).
\]
Consequently,
\begin{align*}
d(x,u)  &  \leq d(x,\overline{x})+d(\overline{x},u)<\tau\frac{\eta}%
{4}+(1+\varepsilon)\frac{\tau}{2}d(z,H(x))\\
&  <\frac{\tau}{2}d(z,H(x))+(1+\varepsilon)\frac{\tau}{2}d(z,H(x)).
\end{align*}
Making $\varepsilon\rightarrow0$, it follows that
\[
d(x,H^{-1}(z))\leq\tau d(z,H(x)),
\]
i.e. the conclusion.$\hfill\square$

\begin{thm}
\label{mainresult} Let $X,Y_{1},Y_{2}$ be complete metric spaces, and $Z$ be a
metric space. Suppose that $F_{1}:X\rightrightarrows Y_{1},$ $F_{2}%
:X\rightrightarrows Y_{2}$ and $G:Y_{1}\times Y_{2}\rightrightarrows Z$ are
closed-graph multifunctions, satisfying the following conditions for some
$(\overline{x},\overline{y}_{1},\overline{y}_{2},\overline{z})\in X\times
Y_{1}\times Y_{2}\times Z$ with $(\overline{x},\overline{y}_{1})\in
\operatorname*{Gr}F_{1},$ $(\overline{x},\overline{y}_{2})\in
\operatorname*{Gr}F_{2},$ $((\overline{y}_{1},\overline{y}_{2}),\overline
{z})\in\operatorname*{Gr}G:$

{(i)} $F_{1}$ is metrically regular around $(\overline{x},\overline{y}_{1})$;

{(ii)} $F_{2}$ has the Aubin property around $(\overline{x},\overline{y}_{2})$;

{(iii)} $G$ is metrically regular around $((\overline{y}_{1},\overline{y}%
_{2}),\overline{z})$ with respect to $y_{1},$ uniformly in $y_{2}$;

{(iv)} $G$ has the Aubin property around $((\overline{y}_{1},\overline{y}%
_{2}),\overline{z})$ with respect to $y_{2},$ uniformly in $y_{1}$;

{(v)} $0<\operatorname*{reg}F_{1}(\overline{x},\overline{y}_{1})\cdot
\operatorname*{lip}F_{2}(\overline{x},\overline{y}_{2})\cdot\widehat
{\operatorname*{reg}}_{y_{1}}G((\overline{y}_{1},\overline{y}_{2}%
),\overline{z})\cdot\widehat{\operatorname*{lip}}_{y_{2}}G((\overline{y}%
_{1},\overline{y}_{2}),\overline{z})<1.$

Denote%
\begin{equation}
\rho_{0}:=\dfrac{\operatorname*{reg}F_{1}(\overline{x},\overline{y}_{1}%
)\cdot\widehat{\operatorname*{reg}}_{y_{1}}G((\overline{y}_{1},\overline
{y}_{2}),\overline{z})}{1-\operatorname*{reg}F_{1}(\overline{x},\overline
{y}_{1})\cdot\operatorname*{lip}F_{2}(\overline{x},\overline{y}_{2}%
)\cdot\widehat{\operatorname*{reg}}_{y_{1}}G((\overline{y}_{1},\overline
{y}_{2}),\overline{z})\cdot\widehat{\operatorname*{lip}}_{y_{2}}%
G((\overline{y}_{1},\overline{y}_{2}),\overline{z})}>0. \label{ro0}%
\end{equation}

Then exist a neighborhood $\mathcal{U}\times\mathcal{V}_{1}\times
\mathcal{V}_{2}\times\mathcal{W}\subset X\times Y_{1}\times Y_{2}\times Z$ of
$(\overline{x},\overline{y}_{1},\overline{y}_{2},\overline{z}),$ $\rho>0$
arbitrarily close to $\rho_{0}$ such that $\rho>\rho_{0},$ and an equivalent
metric $d_{0}$ on $X\times Y_{1}\times Y_{2}$ such that%
\begin{equation}
d_{0}((x,y_{1},y_{2}),R^{-1}(z))\leq\rho\cdot\varphi_{R}((x,y_{1}%
,y_{2}),z)\quad\forall(x,y_{1},y_{2},z)\in\mathcal{U}\times\mathcal{V}%
_{1}\times\mathcal{V}_{2}\times\mathcal{W}. \label{aim}%
\end{equation}

As consequence, $R$ is metrically regular around $((\overline{x},\overline
{y}_{1},\overline{y}_{2}),\overline{z})$ with respect to the metric $d_{0}$ on
$X\times Y_{1}\times Y_{2},$ and%
\[
\operatorname*{reg}R((\overline{x},\overline{y}_{1},\overline{y}%
_{2}),\overline{z})\leq\rho_{0}.
\]

Moreover, if $(F_{1},F_{2}),G$ are locally composition-stable around
$(\overline{x},(\overline{y}_{1},\overline{y}_{2}),\overline{z}),$ then $H$ is
metrically regular around $(\overline{x},\overline{z}),$ and%
\[
\operatorname*{reg}H(\overline{x},\overline{z})\leq\rho_{0}.
\]

\end{thm}

\noindent\textbf{Proof.} Take $l,m,\eta,\lambda>0$ such that, for $\xi>0$
arbitrary small,
\begin{align}
m-\xi &  >\operatorname*{reg}F_{1}(\overline{x},\overline{y}_{1}),\nonumber\\
l-\xi &  >\operatorname*{lip}F_{2}(\overline{x},\overline{y}_{2}),\nonumber\\
\eta-\xi &  >\widehat{\operatorname*{lip}}_{y_{2}}G((\overline{y}%
_{1},\overline{y}_{2}),\overline{z}),\label{const}\\
\lambda &  >\widehat{\operatorname*{reg}}_{y_{1}}G((\overline{y}_{1}%
,\overline{y}_{2}),\overline{z}),\nonumber\\
ml\lambda\gamma &  <1.\nonumber
\end{align}

\noindent Endow the product space $X\times Y_{1}\times Y_{2}$ with the
equivalent metric defined by
\[
d_{0}((x,y_{1},y_{2}),(u,v_{1},v_{2}))=\max\{d(x,u),md(y_{1},v_{1}%
),l^{-1}d(y_{2},v_{2})\},
\]

\noindent and set%
\[
\rho:=\dfrac{m\lambda}{1-ml\lambda\eta}>\rho_{0}.
\]

\noindent Notice that $(X\times Y_{1}\times Y_{2},d_{0})$ is a complete metric space.

\noindent Fix $\xi>0$ arbitrary small such that (\ref{const}) holds. By
$(ii),$ there is $\delta_{1}>0$ such that
\begin{equation}
e(F_{2}(x)\cap B(\overline{y}_{2},\delta_{1}),F_{2}(x^{\prime}))\leq
(l-\xi)d(x,x^{\prime})\quad\forall x,x^{\prime}\in B(\overline{x},\delta_{1}).
\label{AubF2}%
\end{equation}

\noindent By $(iv),$ there is $\delta_{2}>0$ such that
\begin{equation}
e(G(y_{1},y_{2})\cap B(\overline{z},\delta_{2}),G(y_{1},y_{2}^{\prime}%
))\leq(\eta-\xi)d(y_{2},y_{2}^{\prime})\quad\forall y_{1}\in B(\overline
{y}_{1},\delta_{2}),(y_{2},y_{2}^{\prime})\in B(\overline{y}_{2},\delta_{2}).
\label{AubG}%
\end{equation}

\noindent Moreover, by $(iii)$, there is $\delta_{3}>0$ such that
\begin{equation}
d(y_{1},G_{y_{2}}^{-1}(z))\leq\lambda d(z,G(y_{1},y_{2}))\quad\forall y_{1}\in
B(\overline{y}_{1},\delta_{3}),y_{2}\in B(\overline{y}_{2},\delta_{3}),z\in
B(\overline{z}.\delta_{3}). \label{mregG}%
\end{equation}

\noindent Finally, by $(i),$ there is $\delta_{4}>0$ such that
\begin{equation}
d(x,F_{1}^{-1}(y_{1}))\leq(m-\xi)d(y_{1},F_{1}(x))\quad\forall x\in
B(\overline{x},\delta_{4}),y_{1}\in B(\overline{y}_{1},\delta_{4}).
\label{mregF1}%
\end{equation}

\noindent Set
\begin{align*}
\delta &  :=\min\{2^{-1}\delta_{1},2^{-1}\delta_{2},\delta_{3},2^{-1}%
\delta_{4}\},\\
\gamma &  :=\min\left\{  \dfrac{\delta}{\lambda},\dfrac{\delta}{m\lambda
},\dfrac{\delta}{m\lambda l}\right\}  ,
\end{align*}
and take arbitrary $(x,y_{1},y_{2},z)\in B(\overline{x},\delta)\times
B(\overline{y}_{1},\delta)\times B(\overline{y}_{2},\delta)\times
B(\overline{z},\delta),$ such that $z\notin R(x,y_{1},y_{2}),$ $\varphi
_{R}((x,y_{1},y_{2}),z)<\gamma.$ Then $(y_{1},y_{2})\in F_{1}(x)\times
F_{2}(x),$ because otherwise $\varphi_{R}((x,y_{1},y_{2}),z)=+\infty.$ Hence
$z\not \in G(y_{1},y_{2}).$

We want to prove next that for every $\tau>0,$ there exists $(u,v_{1}%
,v_{2})\in X\times Y_{1}\times Y_{2}$ such that%
\begin{equation}
0<d_{0}((x,y_{1},y_{2}),(u,v_{1}v_{2}))<(\rho+\tau)\left[  \varphi
_{R}((x,y_{1},y_{2}),z)-\varphi_{R}((u,v_{1},v_{2}),z)\right]  .
\label{interm1}%
\end{equation}

\noindent For this, fix $\tau>0$ and take arbitrary $\{(x_{n},y_{1n}%
,y_{2n})\}\subset X\times Y_{1}\times Y_{2}$ converging to $(x,y_{1},y_{2})$
such that%
\[
\lim\limits_{n\rightarrow\infty}d(z,R(x_{n},y_{1n},y_{2n}))=\varphi
_{R}((x,y_{1},y_{2}),z)<\gamma.
\]

Hence, as above, for every $n$ sufficiently large, $(y_{1n},y_{2n})\in
F_{1}(x_{n})\times F_{2}(x_{n}).$ Also, using the closedness of
$\operatorname*{Gr}G,$ it follows that, for every $n$ sufficiently large,%
\[
z\not \in G(y_{1n},y_{2n})\text{ and }d(z,G(y_{1n},y_{2n}))<\gamma.
\]

Now, because $(y_{1n},y_{2n})\rightarrow(y_{1},y_{2}),$ one deduces that
$y_{1n}\in B(\overline{y}_{1},\delta)$ and $y_{2n}\in B(\overline{y}%
_{2},\delta),$ for every $n$ sufficiently large. We can use now (\ref{mregG})
to get that%
\begin{equation}
d(y_{1n},G_{y_{2n}}^{-1}(z))\leq\lambda d(z,G(y_{1n},y_{2n}))<\lambda
\gamma\leq\delta. \label{nonempt}%
\end{equation}

Remark that the closedness of $\operatorname*{Gr}G$ assures the fact that
$d(y_{1n},G_{y_{2n}}^{-1}(z))>0$ for every $n\ $large enough. Also, by
(\ref{nonempt}), one may suppose that $G_{y_{2n}}^{-1}(z)\not =\emptyset,$ so
for every $\varepsilon>0,$ there exists $v_{1n}\in G_{y_{2n}}^{-1}(z),$ or
$z\in G(v_{1n},y_{2n}),$ such that%
\begin{align}
d(y_{1n},v_{1n})  &  <\left(  1+\frac{\varepsilon}{2\lambda}\right)
d(y_{1n},G_{y_{2n}}^{-1}(z))\leq\left(  \lambda+\frac{\varepsilon}{2}\right)
d(z,G(y_{1n},y_{2n}))\nonumber\\
&  =\left(  \lambda+\frac{\varepsilon}{2}\right)  [d(z,G(y_{1n},y_{2n}%
))-d(z,G(v_{1n},y_{2n})]. \label{estG}%
\end{align}
\noindent Observe that we may suppose without loosing the generality that
$\lim\limits_{n\rightarrow\infty}d(v_{1n},y_{1n})$ exists and $\lim
\limits_{n\rightarrow\infty}d(v_{1n},y_{1n})=\lim\limits_{n\rightarrow\infty
}d(v_{1n},y_{1})>0,$ because otherwise we would obtain using the closedness of
$\operatorname*{Gr}G$ that $z\in G(y_{1},y_{2}).$

Also, for every $\varepsilon>0$ sufficiently small such that $\left(
1+\dfrac{\varepsilon}{2\lambda}\right)  d(z,G(y_{1n},y_{2n}))<\gamma,$ one has
that $d(y_{1n},v_{1n})<\lambda\gamma\leq\delta,$ so%
\[
d(v_{1n},\overline{y}_{1})\leq d(v_{1n},y_{1n})+d(y_{1n},\overline{y}%
_{1})<2\delta\leq\delta_{4}.
\]

As $x_{n}\in B(\overline{x},\delta_{4}),$ $v_{1n}\in B(\overline{y}_{1}%
,\delta_{4})$ and $d(y_{1n},v_{1n})>0$ for every $n$ sufficiently large, by
(\ref{mregF1}), we have that
\[
d(x_{n},F_{1}^{-1}(v_{1n})\leq(m-\xi)d(v_{1n},F_{1}(x_{n}))<md(v_{1n}%
,y_{1n})<m\lambda\gamma\leq\delta,
\]

\noindent so there is $u_{n}\in F_{1}^{-1}(v_{1n})$ such that%
\begin{equation}
d(x_{n},u_{n})<md(v_{1n},y_{1n})<\delta. \label{est1}%
\end{equation}

\noindent Hence,%
\[
d(u_{n},\overline{x})\leq d(u_{n},x_{n})+d(x_{n},\overline{x})<2\delta
\leq\delta_{1},
\]

\noindent or $u_{n}\in B(\overline{x},\delta_{1}).$ Since $x_{n}\in
B(\overline{x},\delta_{1}),$ $y_{2n}\in F_{2}(x_{n})\cap B(\overline{y}%
_{2},\delta_{1})$, using (\ref{AubF2}), we have that%
\[
d(y_{2n},F_{2}(u_{n}))\leq e(F_{2}(x_{n})\cap B(\overline{y}_{2},\delta
_{1}),F_{2}(u_{n}))\leq(l-\xi)d(x_{n},u_{n}),
\]

\noindent so there exists $v_{2n}\in F_{2}(u_{n})$ such that
\begin{equation}
d(y_{2n},v_{2n})\leq ld(x_{n},u_{n})<lmd(v_{1n},y_{1n})<lm\lambda\gamma
\leq\delta. \label{est2}%
\end{equation}

\noindent Then,%
\[
d(v_{2n},\overline{y}_{2})\leq d(y_{2n},v_{2n})+d(v_{2n},\overline{y}%
_{2})<2\delta\leq\delta_{2}.
\]

\noindent Finally, as $d(v_{1n},\overline{y}_{1})<2\delta\leq\delta_{1},$ by
(\ref{AubG}), one has%

\begin{align}
d(z,G(v_{1n},v_{2n}))  &  \leq d(z,G(v_{1n},y_{2n}))+e(G(v_{1n},y_{2n})\cap
B(\overline{z},\delta_{2}),G(v_{1n},v_{2n}))\nonumber\\
&  \leq d(z,G(v_{1n},y_{2n}))+(\eta-\xi)d(y_{2n},v_{2n}). \label{E2}%
\end{align}

Remark that, by (\ref{est1}) and (\ref{est2}), we have that%
\[
l^{-1}d(y_{2n},v_{2n})\leq d(x_{n},u_{n})\leq md(v_{1n},y_{1n}),
\]

\noindent hence%
\[
d_{0}((x_{n},y_{1n},y_{2n}),(u_{n},v_{1n},v_{2n}))=md(y_{1n},v_{1n}%
)>0\quad\forall n\text{ sufficiently large.}%
\]

Without loosing the generality, we may suppose that $\varepsilon>0$ is chosen
sufficiently small such that%
\[
\frac{1}{m(\lambda+\varepsilon/2)}-l\eta>\frac{1}{\rho+\tau/2}.
\]

Also, because $(y_{1n},y_{2n})\in F_{1}(x_{n})\times F_{2}(x_{n})\ $and
$(v_{1n},v_{2n})\in F_{1}(u_{n})\times F_{2}(u_{n}),$ and using also
(\ref{E2}), (\ref{estG}), one gets that%
\begin{align*}
&  \limsup_{n\rightarrow\infty}\frac{d(z,R(x_{n},y_{1n},y_{2n}))-d(z,R(u_{n}%
,v_{1n},v_{2n}))}{d_{0}((x_{n},y_{1n},y_{2n}),(u_{n},v_{1n},v_{2n}))}\\
&  \geq\limsup_{n\rightarrow\infty}\frac{d(z,G(y_{1n},y_{2n}))-d(z,G(v_{1n}%
,v_{2n}))}{d_{0}((x_{n},y_{1n},y_{2n}),(u_{n},v_{1n},v_{2n}))}\\
&  \geq\limsup_{n\rightarrow\infty}\frac{d(z,G(y_{1n},y_{2n}))-d(z,G(v_{1n}%
,y_{2n}))-\eta d(y_{2n},v_{2n})}{\max\{d(x_{n},u_{n}),md(y_{1n},v_{1n}%
),l^{-1}d(y_{2n},v_{2n})\}}\\
&  \geq\limsup_{n\rightarrow\infty}\frac{d(z,G(y_{1n},y_{2n}))-d(z,G(v_{1n}%
,y_{2n}))}{md(y_{1n},v_{1n})}-\frac{\eta d(y_{2n},v_{2n})}{l^{-1}%
d(y_{2n},v_{2n})}\\
&  \geq\frac{1}{m(\lambda+\varepsilon/2)}-l\eta>\frac{1}{\rho+\tau/2}.
\end{align*}

Remark that
\begin{align*}
\lim_{n\rightarrow\infty}d_{0}((u_{n},v_{1n},v_{2n}),(x,y_{1},y_{2}))  &
=\lim_{n\rightarrow\infty}d_{0}((u_{n},v_{1n},v_{2n}),(x_{n},y_{1n},y_{2n}))\\
&  =\lim_{n\rightarrow\infty}md(y_{1n},v_{1n})>0.
\end{align*}
Take now $r\in\left(  0,\lim_{n\rightarrow\infty}d_{0}((u_{n},v_{1n}%
,v_{2n}),(x,y_{1},y_{2}))\right)  .$ Because $\lim_{\theta\rightarrow0_{+}%
}\dfrac{1+\theta}{1-\theta}\left(  \rho+\dfrac{\tau}{2}\right)  =\rho
+\dfrac{\tau}{2}<\rho+\tau,$ there exists $\theta>0$ sufficiently small such
that, for every $n$ sufficiently large,
\begin{align*}
\dfrac{1+\theta}{1-\theta}\left(  \rho+\dfrac{\tau}{2}\right)   &  <\rho
+\tau,\\
d_{0}((u_{n},v_{1n},v_{2n}),(x_{n},y_{1n},y_{2n}))  &  \geq r,\quad
d_{0}((x_{n},y_{1n},y_{2n}),(x,y_{1},y_{2}))<\theta r,\\
d(z,R(x_{n},y_{1n},y_{2n}))  &  <\varphi_{R}((x,y_{1},y_{2}),z)+\frac{\theta
}{\rho+\dfrac{\tau}{2}}d_{0}((u_{n},v_{1n},v_{2n}),(x_{n},y_{1n},y_{2n})),\\
d_{0}((u_{n},v_{1n},v_{2n}),(x_{n},y_{1n},y_{2n}))  &  <\left(  \rho
+\dfrac{\tau}{2}\right)  \left[  d(z,R(x_{n},y_{1n},y_{2n}))-d(z,R(u_{n}%
,v_{1n},v_{2n}))\right]  .
\end{align*}

It follows that, for every $n$ sufficiently large,%
\begin{align*}
d_{0}((u_{n},v_{1n},v_{2n}),(x_{n},y_{1n},y_{2n}))  &  <\frac{\rho+\dfrac
{\tau}{2}}{1-\theta}\left[  \varphi_{R}((x,y_{1},y_{2}),z)-\varphi_{R}%
((u_{n},v_{1n},v_{2n}),z)\right]  ,\\
d_{0}((u_{n},v_{1n},v_{2n}),(x,y_{1},y_{2}))  &  \leq(1+\theta)d_{0}%
((u_{n},v_{1n},v_{2n}),(x_{n},y_{1n},y_{2n}))\\
&  <\dfrac{1+\theta}{1-\theta}\left(  \rho+\dfrac{\tau}{2}\right)  \left[
\varphi_{R}((x,y_{1},y_{2}),z)-\varphi_{R}((u_{n},v_{1n},v_{2n}),z)\right] \\
&  <(\rho+\tau)\left[  \varphi_{R}((x,y_{1},y_{2}),z)-\varphi_{R}%
((u_{n},v_{1n},v_{2n}),z)\right]  ,
\end{align*}

\noindent i.e. (\ref{interm1}) is proved.

Fix now $t\in(0,\min\{\rho,\tau\}),$ denote%
\[
s:=\min\left\{  \delta,\dfrac{\delta}{6\rho},\dfrac{\delta m}{6\rho}%
,\dfrac{\delta}{6l\rho},\dfrac{\gamma}{4}\right\}  ,
\]
and take arbitrary $z\in B(\overline{z},s)\subset B(\overline{z},\delta).$
Then%
\[
\varphi_{R}((\overline{x},\overline{y}_{1},\overline{y}_{2}),z)\leq
d(z,R(\overline{x},\overline{y}_{1},\overline{y}_{2}))\leq d(z,\overline
{z})<s,
\]

\noindent hence%
\[
\varphi_{R}((\overline{x},\overline{y}_{1},\overline{y}_{2}),z)<\inf
_{(x,y_{1},y_{2})\in X\times Y_{1}\times Y_{2}}\varphi_{R}((x,y_{1}%
,y_{2}),z)+s.
\]

One can apply now the Ekeland Variational Principle to the lower
semicontinuous function $(x,y_{1},y_{2})\mapsto\varphi_{R}((x,y_{1},y_{2}),z)$
and for the distance $d_{0}$ to get the existence of $(u,v_{1},v_{2})\in
X\times Y_{1}\times Y_{2}$ such that%
\begin{align}
\varphi_{R}((u,v_{1},v_{2}),z)  &  \leq\varphi_{R}((\overline{x},\overline
{y}_{1},\overline{y}_{2}),z)<s,\label{ek1}\\
d_{0}((u,v_{1},v_{2}),(\overline{x},\overline{y}_{1},\overline{y}_{2}))  &
\leq s(\rho+t),\label{ek2}\\
\varphi_{R}((u,v_{1},v_{2}),z)  &  \leq\varphi_{R}((x^{\prime},y_{1}^{\prime
},y_{2}^{\prime}),z)+\frac{1}{\rho+t}d_{0}((u,v_{1},v_{2}),(x^{\prime}%
,y_{1}^{\prime},y_{2}^{\prime})), \label{ek3}%
\end{align}

\noindent for every $(x^{\prime},y_{1}^{\prime},y_{2}^{\prime})\in X\times
Y_{1}\times Y_{2}.$

It follows from (\ref{ek2}) that $(u,v_{1},v_{2})\in B(\overline{x}%
,\delta)\times B(\overline{y}_{1},\delta)\times B(\overline{y}_{2},\delta).$
Suppose now that one has $\varphi_{R}((u,v_{1},v_{2}),z)>0.$ By (\ref{ek1}),
$\varphi_{R}((u,v_{1},v_{2}),z)<s<\gamma,$ hence one can deduce the existence
of $(u^{\prime},v_{1}^{\prime},v_{2}^{\prime})\in X\times Y_{1}\times Y_{2}$
such that, using relations (\ref{interm1}) and (\ref{ek3}),%
\begin{align*}
0  &  <d_{0}((u,v_{1},v_{2}),(u^{\prime},v_{1}^{\prime},v_{2}^{\prime}%
))<(\rho+t)\left[  \varphi_{R}((u,v_{1},v_{2}),z)-\varphi_{R}((u^{\prime
},v_{1}^{\prime},v_{2}^{\prime}),z)\right] \\
&  \leq d_{0}((u,v_{1},v_{2}),(u^{\prime},v_{1}^{\prime},v_{2}^{\prime})),
\end{align*}

\noindent a contradiction. Hence, $\varphi_{R}((u,v_{1},v_{2}),z)=0,$ or,
equivalently, $(u,v_{1},v_{2})\in R^{-1}(z).$

Fix now arbitrary
\begin{align*}
(x,y_{1},y_{2})  &  \in B_{0}((\overline{x},\overline{y}_{1},\overline{y}%
_{2}),2s\rho):=\{(x,y_{1},y_{2})\mid d_{0}((x,y_{1},y_{2}),(\overline
{x},\overline{y}_{1},\overline{y}_{2}))<2s\rho\}\\
&  \subset B(\overline{x},\delta)\times B(\overline{y}_{1},\delta)\times
B(\overline{y}_{2},\delta).
\end{align*}
If $\varphi_{R}((x,y_{1},y_{2}),z)\geq\gamma,$ then%
\begin{align*}
d_{0}((x,y_{1},y_{2}),R^{-1}(z))  &  \leq d_{0}((x,y_{1},y_{2}),(u,v_{1}%
,v_{2}))\\
&  \leq d_{0}((x,y_{1},y_{2}),(\overline{x},\overline{y}_{1},\overline{y}%
_{2}))+d_{0}((\overline{x},\overline{y}_{1},\overline{y}_{2}),(u,v_{1}%
,v_{2}))\\
&  <4s\rho\leq\gamma\rho\leq\rho\varphi_{R}((x,y_{1},y_{2}),z).
\end{align*}

\noindent The case $\varphi_{R}((x,y_{1},y_{2}),z)=0$ makes relation
(\ref{aim}) to hold trivially. Consider now the case $\varphi_{R}%
((x,y_{1},y_{2}),z)\in(0,\gamma).$ If there is no $(p,w_{1},w_{2})\in X\times
Y_{1}\times Y_{2}$ such that $d_{0}((x,y_{1},y_{2}),(p,w_{1},w_{2}%
))<d_{0}((x,y_{1},y_{2}),R^{-1}(z))$ and $\varphi_{R}((p,w_{1},w_{2}%
),z)\leq\varphi_{R}((x,y_{1},y_{2}),z),$ then%
\begin{align*}
m(x,y_{1},y_{2})  &  :=\inf\left\{
\begin{array}
[c]{c}%
\sup_{\substack{(u^{\prime},v_{1}^{\prime},v_{2}^{\prime})\in X\times
Y_{1}\times Y_{2},\\(u^{\prime},v_{1}^{\prime},v_{2}^{\prime})\not =%
(x,y_{1},y_{2})}}\dfrac{\varphi_{R}((p,w_{1},w_{2}),z)-\varphi_{R}((u^{\prime
},v_{1}^{\prime},v_{2}^{\prime}),z)}{d_{0}((p,w_{1},w_{2}),(u^{\prime}%
,v_{1}^{\prime},v_{2}^{\prime}))}\mid\\%
\begin{array}
[c]{c}%
d_{0}((x,y_{1},y_{2}),(p,w_{1},w_{2}))<d_{0}((x,y_{1},y_{2}),R^{-1}(z))\\
\varphi_{R}((p,w_{1},w_{2}),z)\leq\varphi_{R}((x,y_{1},y_{2}),z)
\end{array}
\end{array}
\right\} \\
&  =+\infty>\frac{1}{\rho+\tau}.
\end{align*}

Take now any $(p,w_{1},w_{2})\in X\times Y_{1}\times Y_{2}$ such that
$d_{0}((x,y_{1},y_{2}),(p,w_{1},w_{2}))<d_{0}((x,y_{1},y_{2}),R^{-1}(z))$ and
$\varphi_{R}((p,w_{1},w_{2}),z)\leq\varphi_{R}((x,y_{1},y_{2}),z).$ Then%
\begin{align*}
d_{0}((p,w_{1},w_{2}),(\overline{x},\overline{y}_{1},\overline{y}_{2}))  &
\leq d_{0}((p,w_{1},w_{2}),(x,y_{1},y_{2}))+d_{0}((x,y_{1},y_{2}%
),(\overline{x},\overline{y}_{1},\overline{y}_{2}))\\
&  <d_{0}((x,y_{1},y_{2}),R^{-1}(z))+d_{0}((x,y_{1},y_{2}),(\overline
{x},\overline{y}_{1},\overline{y}_{2}))\\
&  \leq d_{0}((\overline{x},\overline{y}_{1},\overline{y}_{2}),R^{-1}%
(z))+2d_{0}((x,y_{1},y_{2}),(\overline{x},\overline{y}_{1},\overline{y}%
_{2}))\\
&  \leq d_{0}((u,v_{1},v_{2}),(\overline{x},\overline{y}_{1},\overline{y}%
_{2}))+2d_{0}((x,y_{1},y_{2}),(\overline{x},\overline{y}_{1},\overline{y}%
_{2}))\leq6\rho s,
\end{align*}

\noindent hence $(p,w_{1},w_{2})\in B_{0}((\overline{x},\overline{y}%
_{1},\overline{y}_{2}),6\rho s)\subset B(\overline{x},\delta)\times
B(\overline{y}_{1},\delta)\times B(\overline{y}_{2},\delta).$ According to
(\ref{interm1}), there exists $(u^{\prime\prime},v_{1}^{\prime\prime}%
,v_{2}^{\prime\prime})\in X\times Y_{1}\times Y_{2},$ $(u^{\prime\prime}%
,v_{1}^{\prime\prime},v_{2}^{\prime\prime})\not =(p,w_{1},w_{2}),$ such that%
\[
\frac{\varphi_{R}((p,w_{1},w_{2}),z)-\varphi_{R}((u^{\prime\prime}%
,v_{1}^{\prime\prime},v_{2}^{\prime\prime}),z)}{d_{0}((p,w_{1},w_{2}%
),(u^{\prime\prime},v_{1}^{\prime\prime},v_{2}^{\prime\prime}))}>\frac{1}%
{\rho+\tau}.
\]

But this means that $m(x,y_{1},y_{2})\geq\dfrac{1}{\rho+\tau}.$ As
$\varphi_{R}((\cdot,\cdot,\cdot),z)$ is lower semicontinuous and
$\{(a,b_{1},b_{2})\in X\times Y_{1}\times Y_{2}:\varphi_{R}((a,b_{1}%
,b_{2}),z)\leq0\}=R^{-1}(z)$, apply Theorem \ref{m} to get that%
\[
\dfrac{1}{\rho+\tau}\cdot d_{0}((x,y_{1},y_{2}),R^{-1}(z))\leq\varphi
_{R}((x,y_{1},y_{2}),z).
\]

Take $\mathcal{U}_{1}\times\mathcal{V}_{1}\times\mathcal{V}_{2}\subset
B_{0}((\overline{x},\overline{y}_{1},\overline{y}_{2}),2s\rho)$ a neighborhood
of $(\overline{x},\overline{y}_{1},\overline{y}_{2}).$ As $\tau$ can be made
arbitrary small, $l,m,\eta,\lambda$ can be taken arbitrary close to the
corresponding regularity moduli, and $(x,y_{1},y_{2})\in\mathcal{U}_{1}%
\times\mathcal{V}_{1}\times\mathcal{V}_{2},$ $z\in\mathcal{W}:=B(\overline
{z},s)$ can be taken arbitrarily, we get the first conclusion.

For the second conclusion, take into account that $\varphi_{R}((x,y_{1}%
,y_{2}),z)\leq d(z,R(x,y_{1},y_{2}))$ for every $(x,y_{1},y_{2},z)\in X\times
Y_{1}\times Y_{2}\times Z.$

Now, observe that the implication $(i)\Rightarrow(ii)$ of Proposition
\ref{proMT} can be identically proven if one replaces the additive distance
with $d_{0}.$ In fact, it is essential to be able to write again relation
(\ref{ineg_dist2}) from (\ref{ineg_dist1}) written for $d_{0},$ and this can
be obviously done. Using Proposition \ref{MTcomsta}, one gets the final
conclusion.$\hfill\square$

\begin{rmk}
This Theorem and Proposition \ref{usc-compositionnew} imply directly Theorem
4.6 given by Ioffe \cite{SiamReFi} for metric regularity of the composition
multifunction of the form $G(x,F(x))$, where $F:X\rightrightarrows Y$ is a
set-valued mapping and $G:X\times Y\rightarrow Y$ is a single-valued mapping.
\end{rmk}

Consider now the case $Y_{1}:=Y_{2}:=Y,$ where $Y$ is a normed vector space,
$G(y_{1},y_{2}):=d(y_{1},y_{2}).$ If $\overline{y}_{1}\not =\overline{y}_{2},$
then $G$ satisfies the conditions from the previous theorem (see \cite[Section
4]{DurStr2012}). Observe now that%
\begin{align*}
H^{-1}(0)  &  =\left\{  x\in X\mid0\in%
{\displaystyle\bigcup\limits_{(y_{1},y_{2})\in F_{1}(x)\times F_{2}(x)}}
d(y_{1},y_{2})\right\} \\
&  =\left\{  x\in X\mid\exists y\in Y\text{ such that }y\in F_{1}(x)\text{ and
}y\in F_{2}(x)\right\} \\
&  =\operatorname*{Fix}(F_{1}^{-1}\circ F_{2}).
\end{align*}

Also,%
\begin{align*}
d(0,G(F_{1}(x)\cap\mathcal{V}_{1},F_{2}(x)\cap\mathcal{V}_{2}))  &  =d\left(
0,%
{\displaystyle\bigcup\limits_{(y_{1},y_{2})\in F_{1}(x)\cap\mathcal{V}%
_{1}\times F_{2}(x)\cap\mathcal{V}_{2}}}
d(y_{1},y_{2})\right) \\
&  =\inf_{(y_{1},y_{2})\in F_{1}(x)\cap\mathcal{V}_{1}\times F_{2}%
(x)\cap\mathcal{V}_{2}}d(y_{1},y_{2})\\
&  =d\left(  F_{1}(x)\cap\mathcal{V}_{1},F_{2}(x)\cap\mathcal{V}_{2}\right)  .
\end{align*}

As consequence, relation (\ref{aim}) from Theorem \ref{mainresult} implies, in
the virtue of $(i)\Rightarrow(ii)$ from Proposition \ref{proMT}, that there
exist a neighborhood \ $\mathcal{U}\times\mathcal{V}_{1}\times\mathcal{V}%
_{2}\subset X\times Y\times Y$ of $(\overline{x},\overline{y}_{1},\overline
{y}_{2})$ and $\rho>0$ arbitrarily close to $\rho_{0}$ such that
\[
d(x,\operatorname*{Fix}(F_{1}^{-1}\circ F_{2}))\leq\rho d\left(  F_{1}%
(x)\cap\mathcal{V}_{1},F_{2}(x)\cap\mathcal{V}_{2}\right)  \quad\text{for all
}x\in\mathcal{U}.
\]

In this way, one can obtain a direct proof based on Ekeland Variational
Principle, for fixed points/coincidence results of the type Arutyunov, Ioffe
and Dontchev and Frankovska proved recently in \cite{Arut2007},
\cite{Arut2009}, \cite{Ioffe2010b}, \cite{DonFra2010}, \cite{DonFra2}.

\section{Coderivative conditions}

In this section, unless otherwise stated, we assume that all the spaces
involved are Asplund, i.e., Banach spaces where every convex continuous
function is generically Fr\'{e}chet differentiable (in particular, any
reflexive space is Asplund; see, e.g., \cite{Mor2006} for more details). We
recall next some standard notations and definitions we use in the sequel. If
$X$ is a normed vector space, we denote by $B_{X},$ $\overline{B}_{X},$
$S_{X}$ the open unit ball, the closed unit ball and the unit sphere of $X,$
respectively. As usual, by $X^{\ast}$ we denote the topological dual of the
normed vector space $X,$ while the symbol $\left\langle \cdot,\cdot
\right\rangle $ stands for the canonical duality pairing between $X$ and
$X^{\ast}$. The symbol $\overset{w^{\ast}}{\rightarrow}$ indicates the
convergence in the weak-star topology of $X^{\ast}.$ Given a set-valued
mapping $F\colon X\rightrightarrows X^{\ast}$, recall that
\begin{equation}
\underset{x\rightarrow\overline{x}}{\operatorname{Limsup}}F(x):=\left\{
x^{\ast}\in X^{\ast}\mid\exists x_{n}\rightarrow\overline{x},x_{n}^{\ast
}\overset{w^{\ast}}{\rightarrow}x^{\ast}\;\text{with}\;x_{n}^{\ast}\in
F(x_{n}),\text{ }n\in\mathbb{N}\right\}  \label{PK}%
\end{equation}
stands for the sequential Painlev\'{e}-Kuratowski outer/upper limit of $F$ as
$x\rightarrow\overline{x}$ with respect to the norm topology of $X$ and the
weak$^{\ast}$ topology of $X^{\ast}$.

If $X$ is a normed vector space, $S$ is a non-empty subset of $X$ and $x\in
S,$ $\varepsilon\geq0,$ the set of $\varepsilon-$normals to $S$ at $x$ is%
\begin{equation}
\widehat{N}_{\varepsilon}(S,x):=\left\{  x^{\ast}\in X^{\ast}\mid
\underset{u\overset{S}{\rightarrow}x}{\lim\sup}\frac{\left\langle u-x,x^{\ast
}\right\rangle }{\left\Vert u-x\right\Vert }\leq\varepsilon\right\}  .
\label{eps-no}%
\end{equation}

\noindent If $\varepsilon=0,$ the elements in the right-hand side of
(\ref{eps-no}) are called Fr\'{e}chet normals and their collection, denoted by
$\widehat{N}(S,x),$ is the Fr\'{e}chet normal cone to $S$ at $x.$

Let $\overline{x}\in S.$ The basic (or limiting, or Mordukhovich) normal cone
to $S$ at $\overline{x}$ is
\[
N(S,\overline{x}):=\{x^{\ast}\in X^{\ast}\mid\exists\varepsilon_{n}%
\downarrow0,x_{n}\overset{S}{\rightarrow}\overline{x},x_{n}^{\ast}%
\overset{w^{\ast}}{\rightarrow}x^{\ast},x_{n}^{\ast}\in\widehat{N}%
_{\varepsilon_{n}}(S,x_{n}),\forall n\in%
\mathbb{N}
\},
\]

\noindent where by $x_{n}\overset{S}{\rightarrow}\overline{x}$ we denote
$x_{n}\overset{S}{\rightarrow}\overline{x},x_{n}\in S$ for every $n$
sufficiently large.

If $X$ is an Asplund space and $S$ is locally closed around $\overline{x}$,
the formula for the basic normal cone takes a simpler form:%
\[
N(S,\overline{x})=\{x^{\ast}\in X^{\ast}\mid\exists x_{n}\overset
{S}{\rightarrow}\overline{x},x_{n}^{\ast}\overset{w^{\ast}}{\rightarrow
}x^{\ast},x_{n}^{\ast}\in\widehat{N}(S,x_{n}),\forall n\in%
\mathbb{N}
\},
\]

\noindent i.e.%
\[
N(S,\overline{x})=\underset{x\overset{S}{\rightarrow}\overline{x}%
}{\operatorname{Limsup}}\widehat{N}(S,x).
\]

Consider now $f:X\rightarrow\overline{\mathbb{R}}\ $such that is finite at
$\overline{x}\in X.$ The Fr\'{e}chet subdifferential of $f$ at $\overline{x}$
is the set
\[
\widehat{\partial}f(\overline{x}):=\{x^{\ast}\in X^{\ast}\mid(x^{\ast}%
,-1)\in\widehat{N}(\operatorname*{epi}f,(\overline{x},f(\overline{x})))\},
\]
and the basic (or limiting, or Mordukhovich) subdifferential of $f$ at
$\overline{x}$ is%
\[
\partial f(\overline{x}):=\{x^{\ast}\in X^{\ast}\mid(x^{\ast},-1)\in
N(\operatorname*{epi}f,(\overline{x},f(\overline{x})))\},
\]
where $\operatorname*{epi}f$ denotes the epigraph of $f.$

In general Banach spaces, for every lower semicontinuous function one has
$\widehat{\partial}f(\overline{x})\subset\partial f(\overline{x}),$ and in
Asplund spaces the next relation holds
\[
\partial f(\overline{x})=\underset{x\overset{f}{\rightarrow}\overline{x}%
}{\operatorname{Limsup}}\widehat{\partial}f(x),
\]
where by $x\overset{f}{\rightarrow}\overline{x}$ we mean $x\rightarrow
\overline{x},f(x)\rightarrow f(\overline{x}).$

It is well-known that the Fr\'{e}chet subdifferential satisfies a fuzzy sum
rule on Asplund spaces (see \cite[Theorem 2.33]{Mor2006}). More precisely, if
$X$ is an Asplund space and $\varphi_{1},\varphi_{2}:X\rightarrow
\mathbb{R\cup\{\infty\}}$ are such that $\varphi_{1}$ is Lipschitz continuous
around $\overline{x}\in\operatorname*{dom}\varphi_{1}\cap\operatorname*{dom}%
\varphi_{2}$ and $\varphi_{2}$ is lower semicontinuous around $\overline{x},$
then for any $\gamma>0$ one has%
\begin{equation}
\widehat{\partial}(\varphi_{1}+\varphi_{2})(\overline{x})\subset%
{\displaystyle\bigcup}
\{\widehat{\partial}\varphi_{1}(x_{1})+\widehat{\partial}\varphi_{2}%
(x_{2})\mid x_{i}\in\overline{x}+\gamma\overline{B}_{X},\left\vert \varphi
_{i}(x_{i})-\varphi_{i}(\overline{x})\right\vert \leq\gamma,i=1,2\}+\gamma
\overline{B}_{X^{\ast}}. \label{fuz}%
\end{equation}

If $\delta_{\Omega}$ denotes the indicator function associated to a nonempty
set $S\subset X$ (i.e. $\delta_{S}(x)=0$ if $x\in S,$ $\delta_{S}(x)=\infty$
otherwise), then for any $\overline{x}\in S,$ $\widehat{\partial}\delta
_{S}(\overline{x})=\widehat{N}(S,\overline{x})$ and $\partial\delta
_{S}(\overline{x})=N(S,\overline{x}).$

Finally, let $F:X\rightrightarrows Y$ be a set-valued map and $(\overline
{x},\overline{y})\in\operatorname*{Gr}F.$ Then the Fr\'{e}chet coderivative at
$(\overline{x},\overline{y})$ is the set-valued map $\widehat{D}^{\ast
}F(\overline{x},\overline{y}):Y^{\ast}\rightrightarrows X^{\ast}$ given by
\[
\widehat{D}^{\ast}F(\overline{x},\overline{y})(y^{\ast}):=\{x^{\ast}\in
X^{\ast}\mid(x^{\ast},-y^{\ast})\in\widehat{N}(\operatorname{Gr}%
F,(\overline{x},\overline{y}))\}.
\]
Similarly, the normal coderivative of $F$ at $(\overline{x},\overline{y})$ is
the set-valued map $D^{\ast}F(\overline{x},\overline{y}):Y^{\ast
}\rightrightarrows X^{\ast}$ given by
\[
D^{\ast}F(\overline{x},\overline{y})(y^{\ast}):=\{x^{\ast}\in X^{\ast}%
\mid(x^{\ast},-y^{\ast})\in N(\operatorname{Gr}F,(\overline{x},\overline
{y}))\}.
\]

\noindent Note that, in fact, the concept of normal coderivative,
independently of the normal cone used in its definition, was introduced in
\cite{Mor1980}.

\bigskip

In order to obtain a result concerning coderivative conditions which assure
the metric regularity of the composition multifunction, we need a calculus
rule for the Fr\'{e}chet normal cone of the intersection of a finite number of
sets. We introduce next some concepts which are needed in order to formulate
such a result.

Recall that the sets $S_{1},...,S_{k}$ satisfy the metric inequality at
$\overline{x}$ if there are $\tau>0$ and $r>0$ such that
\[
d(x,S_{1}\cap...\cap S_{k})\leq\tau\lbrack d(x,S_{1})+...+d(x,S_{k}%
)]\;\text{for every}\;x\in B(\overline{x},r).
\]

\noindent This inequality is shown to be very effective as assumption in order
to infer generalized differentiation calculus rules for different kinds of
operations with sets and multifunctions (see, e.g., \cite{Ioffe2000},
\cite{Pen1998}, \cite{NT2001}).

Given the closed subsets $S_{1},...,S_{k}$ of a normed vector space $X,$ one
says that they are allied at $\overline{x}\in S_{1}\cap...\cap S_{k}$ (for the
Fr\'{e}chet normal cones) whenever $(x_{in})\overset{S_{i}}{\rightarrow
}\overline{x},x_{in}^{\ast}\in\widehat{N}(S_{i},x_{in}),i=\overline{1,k},$ the
relation $(x_{1n}^{\ast}+...+x_{kn}^{\ast})\rightarrow0$ implies
$(x_{in}^{\ast})\rightarrow0$ for every $i=\overline{1,k}$ (for more details,
see \cite{Pen1998}, \cite{LPX2011} and the references therein).

\begin{thm}
Suppose $X$ is an Asplund space and let $S_{1},...,S_{k}$ be closed subsets
such that $\overline{x}\in S_{1}\cap...\cap S_{k}.$ Consider the assertions:

(i) $S_{2},...,S_{k}$ are sequentially normally compact and%
\[
x_{i}^{\ast}\in N(S_{i},\overline{x}),i=\overline{1,k}\text{ and }x_{1}^{\ast
}+...+x_{k}^{\ast}=0\implies x_{1}^{\ast}=...=x_{k}^{\ast}=0;
\]

(ii) $S_{1},...,S_{k}$ are allied at $\overline{x};$

(iii) there exist $a,r>0$ such that for any $x_{i}\in S_{i}\cap B(\overline
{x},r),x_{i}^{\ast}\in\widehat{N}(S_{i},x_{i}),$ one has%
\[
\max_{i=\overline{1,k}}\left\Vert x_{i}^{\ast}\right\Vert \geq1\implies
\left\Vert \sum_{i=1}^{k}x_{i}^{\ast}\right\Vert \geq a^{-1};
\]

(iv) $S_{1},...,S_{k}$ satisfy the metric inequality at $\overline{x};$

(v) there exists $r>0$ such that, for every $\varepsilon>0$ and every
$x\in\lbrack S_{1}\cap...\cap S_{k}]\cap B(\overline{x},r),$ there exist
$x_{i}\in S_{i}\cap B(x,\varepsilon),i=\overline{1,k}$ such that%
\[
\widehat{N}\left(  S_{1}\cap...\cap S_{k},x\right)  \subset\widehat{N}%
(S_{1},x_{1})+...+\widehat{N}(S_{k},x_{k})+\varepsilon\overline{B}_{X^{\ast}%
}.
\]

Then $(i)\Rightarrow(iv)$ and $(ii)\Rightarrow(iii)\Rightarrow(iv)\Rightarrow
(v).$
\end{thm}

\noindent\textbf{Proof.} For the proof of $(i)\Rightarrow(iv),$ see, e.g.,
\cite[Theorem 3.8]{NT2001}, \cite[Corollary 3.10]{Pen1998}. For the rest of
the implications, one may find the proofs (with obvious minor modifications)
in \cite[Theorem 3.7, Propositions 3.8, 3.9]{Pen1998}, \cite[Theorem
3.8]{NT2001}.$\hfill\square$

\bigskip

Let $X,Y_{1},Y_{2},Z$ be Asplund spaces. Suppose that $F_{1}%
:X\rightrightarrows Y_{1},$ $F_{2}:X\rightrightarrows Y_{2}$ and
$G:Y_{1}\times Y_{2}\rightrightarrows Z$ are multifunctions and $(\overline
{x},\overline{y}_{1},\overline{y}_{2},\overline{z})\in X\times Y_{1}\times
Y_{2}\times Z$ is such that $(\overline{x},\overline{y}_{1})\in
\operatorname*{Gr}F_{1},$ $(\overline{x},\overline{y}_{2})\in
\operatorname*{Gr}F_{2},$ $((\overline{y}_{1},\overline{y}_{2}),\overline
{z})\in\operatorname*{Gr}G.$ Setting
\begin{align}
C_{1}  &  :=\{(x,y_{1},y_{2},z)\in X\times Y_{1}\times Y_{2}\times Z:y_{1}\in
F_{1}(x)\},\nonumber\\
C_{2}  &  :=\{(x,y_{1},y_{2},z)\in X\times Y_{1}\times Y_{2}\times Z:y_{2}\in
F_{2}(x)\},\label{C13}\\
C_{3}  &  :=\{(x,y_{1},y_{2},z)\in X\times Y_{1}\times Y_{2}\times Z:z\in
G(y_{1},y_{2})\},\nonumber
\end{align}

\noindent observe that $(\overline{x},\overline{y}_{1},\overline{y}%
_{2},\overline{z})\in C_{1}\cap C_{2}\cap C_{3}.$

Consider now the alliedness property of $C_{1},C_{2},C_{3}$ at $(\overline
{x},\overline{y}_{1},\overline{y}_{2},\overline{z}):$ for any sequences
$(x_{n},y_{1n})\overset{\operatorname*{Gr}F_{1}}{\longrightarrow}(\overline
{x},\overline{y}_{1}),$ $(u_{n},v_{2n})\overset{\operatorname*{Gr}F_{2}%
}{\longrightarrow}(\overline{x},\overline{y}_{2}),$ $(b_{1n},b_{2n}%
,c_{n})\overset{\operatorname*{Gr}G}{\longrightarrow}(\overline{y}%
_{1},\overline{y}_{2},\overline{z})$ and every $x_{n}^{\ast}\in\widehat
{D}^{\ast}F_{1}(x_{n},y_{1n})(y_{1n}^{\ast}),$ $u_{n}^{\ast}\in\widehat
{D}^{\ast}F_{2}(u_{n},v_{2n})(v_{2n}^{\ast}),$ $(b_{1n}^{\ast},b_{2n}^{\ast
})\in\widehat{D}^{\ast}G(b_{1n},b_{2n},c_{n})(c_{n}^{\ast}),$ the relations
$(x_{n}^{\ast}+u_{n}^{\ast})\rightarrow0,$ $\left(  y_{1n}^{\ast}+b_{1n}%
^{\ast}\right)  \rightarrow0,\left(  v_{2n}^{\ast}+b_{2n}^{\ast}\right)
\rightarrow0,\left(  c_{n}^{\ast}\right)  \rightarrow0$ imply
\[
(x_{n}^{\ast})\rightarrow0,(u_{n}^{\ast})\rightarrow0,\left(  y_{1n}^{\ast
})\rightarrow0,(b_{1n}^{\ast}\right)  \rightarrow0,\left(  v_{2n}^{\ast
})\rightarrow0,(b_{2n}^{\ast}\right)  \rightarrow0.
\]

\bigskip

The next result is twofold. On one hand, it provides a sufficient Fr\'{e}chet
coderivative condition for linear openness/metric regularity on Asplund
spaces. On the other hand, it serves in the sequel as the basis for getting
again the conclusion of the main result from the previous section, i.e.
Theorem \ref{mainresult}.

\begin{thm}
\label{coder_Fr} Let $X,Y_{1},Y_{2}$ and $Z$ be Asplund spaces, $F_{1}%
:X\rightrightarrows Y_{1},$ $F_{2}:X\rightrightarrows Y_{2},$ $G:Y_{1}\times
Y_{2}\rightrightarrows Z$ be closed-graph multifunctions, and $(\overline
{x},\overline{y}_{1},\overline{y}_{2},\overline{z})\in X\times Y_{1}\times
Y_{2}\times Z$ be such that $\overline{z}\in G(\overline{y}_{1},\overline
{y}_{2}),(\overline{y}_{1},\overline{y}_{2})\in F_{1}(\overline{x})\times
F_{2}(\overline{x}).$ Assume that the sets $C_{1},C_{2},C_{3}$ defined by
(\ref{C13}) are allied at $(\overline{x},\overline{y}_{1},\overline{y}%
_{2},\overline{z})$ and that there exists $c>0$ such that%
\begin{equation}
c<\liminf_{\substack{{(u_{1},v_{1})\overset{F_{1}}{\rightarrow}(\overline
{x},\overline{y}_{1}),}\text{ }{(u_{2},v_{2})\overset{F_{2}}{\rightarrow
}(\overline{x},\overline{y}_{2})}\\(t_{1},t_{2},w){\overset{G}{\rightarrow
}(\overline{y}}_{1},\overline{y}_{2},\overline{z}{),}\text{ }\delta
\downarrow0}}\left\{  \Vert x_{1}^{\ast}+x_{2}^{\ast}\Vert:\;%
\begin{array}
[c]{l}%
x_{1}^{\ast}\in\widehat{D}^{\ast}F_{1}(u_{1},v_{1})(t_{1}^{\ast}),\\
x_{2}^{\ast}\in\widehat{D}^{\ast}F_{2}(u_{2},v_{2})(t_{2}^{\ast}),\\
(z_{1}^{\ast}+t_{1}^{\ast},z_{2}^{\ast}+t_{2}^{\ast})\in\widehat{D}^{\ast
}G(t_{1},t_{2},w)(w^{\ast}),\\
\Vert w^{\ast}\Vert=1,\Vert z_{1}^{\ast}\Vert<\delta,\Vert z_{2}^{\ast}%
\Vert<\delta,
\end{array}
\right\}  . \tag*{$(C)$}\label{C}%
\end{equation}

Then for every $a\in(0,c),$ $H$ is $a-$open at $(\overline{x},\overline{z}).$
If, moreover, $(F_{1},F_{2}),G$ are locally composition-stable around
$(\overline{x},(\overline{y}_{1},\overline{y}_{2}),\overline{z}),$ then for
every $a\in(0,c),$ $H$ is $a-$open around $(\overline{x},\overline{z}).$
\end{thm}

\noindent\textbf{Proof.} Fix $a\in(0,c).$ Remark that the condition \ref{C}
assures the existence of $r>0$ such that for every $\delta\in(0,r),$ every
${(u_{1},v_{1})\in}\operatorname*{Gr}F_{1}\cap\lbrack B(\overline{x},r)\times
B(\overline{y}_{1},r)],$ ${(u_{2},v_{2})\in}\operatorname*{Gr}F_{2}\cap\lbrack
B(\overline{x},r)\times B(\overline{y}_{2},r)],$ $(t_{1},t_{2},w)\in
\operatorname*{Gr}G\cap\lbrack B(\overline{y}_{1},r)\times B(\overline{y}%
_{2},r)\times B(\overline{z},r)],$ every $w^{\ast}\in S_{Z^{\ast}},$
$t_{1}^{\ast}\in Y_{1}^{\ast},$ $t_{2}^{\ast}\in Y_{2}^{\ast},$ $z_{1}^{\ast
}\in\delta B_{Y_{1}^{\ast}},$ $z_{2}^{\ast}\in\delta B_{Y_{2}^{\ast}}$ such
that $(z_{1}^{\ast}+t_{1}^{\ast},z_{2}^{\ast}+t_{2}^{\ast})\in\widehat
{D}^{\ast}G((t_{1},t_{2}),w)(w^{\ast})$ and every $x_{1}^{\ast}\in\widehat
{D}^{\ast}F_{1}(u_{1},v_{1})(t_{1}^{\ast}),$ $x_{2}^{\ast}\in\widehat{D}%
^{\ast}F_{2}(u_{2},v_{2})(t_{2}^{\ast}),$%
\begin{equation}
c\leq\Vert x_{1}^{\ast}+x_{2}^{\ast}\Vert. \label{ineq_x12}%
\end{equation}

We will prove the existence of $\varepsilon>0$ such that, for every $\rho
\in(0,\varepsilon)$ and every $(\widehat{x},\widehat{y}_{1},\widehat{y}%
_{2},\widehat{z})\in\lbrack B(\overline{x},\varepsilon)\times B(\overline
{y}_{1},\varepsilon)\times B(\overline{y}_{2},\varepsilon)\times
B(\overline{z},\varepsilon)]\cap\lbrack C_{1}\cap C_{2}\cap C_{3}],$%
\begin{equation}
B(\widehat{z},a\rho)\subset H(B(\widehat{x},\rho)). \label{interm}%
\end{equation}

This will show in particular that $H$ is $a-$open at $(\overline{x}%
,\overline{z}).$ Moreover, if $(F_{1},F_{2}),G$ are locally composition-stable
around $(\overline{x},(\overline{y}_{1},\overline{y}_{2}),\overline{z}),$ the
definition of this concept ensures that if $(x,z)$ is close to $(\overline
{x},\overline{z}),$ one can find $(y_{1},y_{2})$ close to $(\overline{y}%
_{1},\overline{y}_{2})$ such that $(x,y_{1},y_{2},z)\in C_{1}\cap C_{2}\cap
C_{3},$ hence one may apply (\ref{interm}) in order to conclude that $H$ is
$a-$open around $(\overline{x},\overline{z}).$

Let us proceed to the proof of (\ref{interm}). Choose $b\in(0,1)$ such that
$\dfrac{a}{a+1}<b<\dfrac{c}{c+1}\ $and fix $\delta\in(0,r)$ as above. Then one
can find $\alpha,\varepsilon>0$ such that the following are satisfied:

\begin{itemize}
\item $b^{-1}a\varepsilon<4^{-1}r;$

\item $\alpha^{-1}a\varepsilon<4^{-1}r;$

\item $\varepsilon<2^{-1}r;$

\item $\dfrac{a}{a+1}<b+2\varepsilon<\dfrac{c}{c+1};$

\item $b+2\varepsilon<1;$

\item $\dfrac{\alpha(1+2\varepsilon)}{1-(b+2\varepsilon)}<\delta.$
\end{itemize}

Fix $\rho\in(0,\varepsilon)$ and $(\widehat{x},\widehat{y}_{1},\widehat{y}%
_{2},\widehat{z})\in B(\overline{x},\varepsilon)\times B(\overline{y}%
_{1},\varepsilon)\times B(\overline{y}_{2},\varepsilon)\times B(\overline
{z},\varepsilon)]\cap\lbrack C_{1}\cap C_{2}\cap C_{3}].$ Choose now $v\in
B(\widehat{z},\rho a).$ Observe that the set $C_{1}\cap C_{2}\cap C_{3}$ is
closed because of the closedness of the graphs of $F_{1},F_{2}$ and $G.$ Endow
the space $X\times Y_{1}\times Y_{2}\times Z$ with the distance%
\[
d((x,y_{1},y_{2},z),(x^{\prime},y_{1}^{\prime},y_{2}^{\prime},z^{\prime
})):=b\left\Vert x-x^{\prime}\right\Vert +\alpha\left\Vert y_{1}-y_{1}%
^{\prime}\right\Vert +\alpha\left\Vert y_{2}-y_{2}^{\prime}\right\Vert
+b\left\Vert z-z^{\prime}\right\Vert
\]
and apply the Ekeland Variational Principle for the function $f:C_{1}\cap
C_{2}\cap C_{3}\rightarrow\mathbb{R},$ $f(x,y_{1},y_{2},z):=\left\Vert
v-z\right\Vert $ and $(\widehat{x},\widehat{y}_{1},\widehat{y}_{2},\widehat
{z})\in\operatorname*{dom}f$\ in order to find $(x_{0},y_{10},y_{20},z_{0})\in
C_{1}\cap C_{2}\cap C_{3}$ such that%
\begin{equation}
\left\Vert v-z_{0}\right\Vert \leq\left\Vert v-\widehat{z}\right\Vert
-(b\left\Vert \widehat{x}-x_{0}\right\Vert +\alpha\left\Vert \widehat{y}%
_{1}-y_{10}\right\Vert +\alpha\left\Vert \widehat{y}_{2}-y_{20}\right\Vert
+b\left\Vert \widehat{z}-z_{0}\right\Vert ) \label{ek_n1}%
\end{equation}
and%
\begin{equation}
\left\Vert v-z_{0}\right\Vert \leq\left\Vert v-z\right\Vert +(b\left\Vert
x-x_{0}\right\Vert +\alpha\left\Vert y_{1}-y_{10}\right\Vert +\alpha\left\Vert
y_{2}-y_{20}\right\Vert +b\left\Vert z-z_{0}\right\Vert ), \label{ek_n2}%
\end{equation}

\noindent for every $(x,y_{1},y_{2},z)\in C_{1}\cap C_{2}\cap C_{3}.$

We have%
\begin{align*}
\left\Vert \widehat{x}-x_{0}\right\Vert +\left\Vert \widehat{z}-z_{0}%
\right\Vert  &  \leq b^{-1}\left\Vert v-\widehat{z}\right\Vert <b^{-1}%
a\rho\leq b^{-1}a\varepsilon<4^{-1}r,\\
\left\Vert \widehat{y}_{1}-y_{10}\right\Vert +\left\Vert \widehat{y}%
_{2}-y_{20}\right\Vert  &  \leq\alpha^{-1}\left\Vert v-\widehat{z}\right\Vert
<\alpha^{-1}a\rho\leq\alpha^{-1}a\varepsilon<4^{-1}r,
\end{align*}
hence%
\begin{align*}
(x_{0},y_{10},y_{20},z_{0})  &  \in B(\widehat{x},4^{-1}r)\times B(\widehat
{y}_{1},4^{-1}r)\times B(\widehat{y}_{2},4^{-1}r)\times B(\widehat{z}%
,4^{-1}r).\\
&  \subset B(\overline{x},2^{-1}r)\times B(\overline{y}_{1},2^{-1}r)\times
B(\overline{y}_{2},2^{-1}r)\times B(\overline{z},2^{-1}r).
\end{align*}

If $v=z_{0},$ then%
\begin{align*}
b\left\Vert \widehat{x}-x_{0}\right\Vert  &  \leq(1-b)\left\Vert v-\widehat
{z}\right\Vert \\
&  <(1-b)a\rho<b\rho,
\end{align*}
hence $x_{0}\in B(\widehat{x},\rho)$ and $v=z_{0}\in H(x_{0})\subset
H(B(\widehat{x},\rho)),$ which proves the desired conclusion.

We want to prove now that $v=z_{0}$ is the sole possible situation. Suppose
then by contradiction that $v\not =z_{0}$ and consider the function $h:X\times
Y_{1}\times Y_{2}\times Z\rightarrow\mathbb{R},$%
\[
h(x,y_{1},y_{2},z):=\left\Vert v-z\right\Vert +(b\left\Vert x-x_{0}\right\Vert
+\alpha\left\Vert y_{1}-y_{10}\right\Vert +\alpha\left\Vert y_{2}%
-y_{20}\right\Vert +b\left\Vert z-z_{0}\right\Vert ).
\]
From (\ref{ek_n2}), we have that the point $(x_{0},y_{10},y_{20},z_{0})$ is a
minimum point for $h$ on the set $C_{1}\cap C_{2}\cap C_{3}$, or,
equivalently, $(x_{0},y_{10},y_{20},z_{0})$ is a global minimum point for the
function $h+\delta_{C_{1}\cap C_{2}\cap C_{3}}.$ Applying the generalized
Fermat rule, we have that%
\[
(0,0,0,0)\in\widehat{\partial}(h(\cdot,\cdot,\cdot,\cdot)+\delta_{C_{1}\cap
C_{2}\cap C_{3}}(\cdot,\cdot,\cdot,\cdot))(x_{0},y_{10},y_{20},z_{0}).
\]

Observe now that $h$ is Lipschitz and $\delta_{C_{1}\cap C_{2}\cap C_{3}}$ is
lower semicontinuous, so one can apply the fuzzy calculus rule for the
Fr\'{e}chet subdifferential. Choose $\gamma\in(0,\min\{\rho,\alpha\rho
,4^{-1}r\})$ such that $v\notin\overline{B}(z_{0},\gamma)$ and obtain that
there exist
\begin{align*}
(u_{1},v_{11},v_{21},w_{1})  &  \in\overline{B}(x_{0},\gamma)\times
\overline{B}(y_{10},\gamma)\times\overline{B}(y_{20},\gamma)\times\overline
{B}(z_{0},\gamma),\\
(u_{2},v_{12},v_{22},w_{2})  &  \in\lbrack\overline{B}(x_{0},\gamma
)\times\overline{B}(y_{10},\gamma)\times\overline{B}(y_{20},\gamma
)\times\overline{B}(z_{0},\gamma)]\cap\lbrack C_{1}\cap C_{2}\cap C_{3}]
\end{align*}
such that%
\begin{equation}
(0,0,0,0)\in\widehat{\partial}h(u_{1},v_{11},v_{21},w_{1})+\widehat{\partial
}\delta_{C_{1}\cap C_{2}\cap C_{3}}(u_{2},v_{12},v_{22},w_{2})+\gamma
(\overline{B}_{X^{\ast}}\times\overline{B}_{Y_{1}^{\ast}}\times\overline
{B}_{Y_{2}^{\ast}}\times\overline{B}_{Z^{\ast}}). \label{fr_sum}%
\end{equation}

Observe that $h$ is the sum of five convex functions, Lipschitz on $X\times
Y_{1}\times Y_{2}\times Z,$ hence $\widehat{\partial}h$ coincides with the sum
of the convex subdifferentials. Also, because $v\not =w_{1}\in\overline
{B}(z_{0},\gamma),$ we get%
\[
\widehat{\partial}\left\Vert v-\cdot\right\Vert (w_{1})=\{z^{\ast}\mid
z^{\ast}\in S_{Z^{\ast}},\text{ }z^{\ast}(v-w_{1})=\left\Vert v-w_{1}%
\right\Vert \}.
\]

Consequently, we have from (\ref{fr_sum}) that%
\begin{align*}
(0,0,0,0)  &  \in\{0\}\times\{0\}\times\{0\}\times S_{Z^{\ast}}+b\overline
{B}_{X^{\ast}}\times\{0\}\times\{0\}\times\{0\}\\
&  +\{0\}\times\alpha\overline{B}_{Y_{1}^{\ast}}\times\{0\}\times
\{0\}+\{0\}\times\{0\}\times\alpha\overline{B}_{Y_{2}^{\ast}}\times
\{0\}+\{0\}\times\{0\}\times\{0\}\times b\overline{B}_{Z^{\ast}}\\
&  +\widehat{N}(C_{1}\cap C_{2}\cap C_{3},(u_{2},v_{12},v_{22},w_{2}%
))+\rho\overline{B}_{X^{\ast}}\times\alpha\rho\overline{B}_{Y_{1}^{\ast}%
}\times\alpha\rho\overline{B}_{Y_{2}^{\ast}}\times\rho\overline{B}_{Z^{\ast}%
}\\
&  \subset\{0\}\times\{0\}\times\{0\}\times S_{Z^{\ast}}+\widehat{N}(C_{1}\cap
C_{2}\cap C_{3},(u_{2},v_{12},v_{22},w_{2}))\\
&  +(b+\rho)\overline{B}_{X^{\ast}}\times\alpha(1+\rho)\overline{B}%
_{Y_{1}^{\ast}}\times\alpha(1+\rho)\overline{B}_{Y_{2}^{\ast}}\times
(b+\rho)\overline{B}_{Z^{\ast}}.
\end{align*}

Now, use the alliedness of $C_{1},C_{2},C_{3}$ at $(\overline{x},\overline
{y}_{1},\overline{y}_{2},\overline{z})$ to get that%
\begin{align*}
\widehat{N}(C_{1}\cap C_{2}\cap C_{3},(u_{2},v_{12},v_{22},w_{2}))  &
\subset\widehat{N}(C_{1},(u_{3},v_{13},v_{23},w_{3}))+\widehat{N}(C_{2}%
,(u_{4},v_{14},v_{24},w_{4}))\\
&  +\widehat{N}(C_{3},(u_{5},v_{15},v_{25},w_{5}))+\rho\overline{B}_{X^{\ast}%
}\times\alpha\rho\overline{B}_{Y_{1}^{\ast}}\times\alpha\rho\overline
{B}_{Y_{2}^{\ast}}\times\rho\overline{B}_{Z^{\ast}},
\end{align*}

\noindent where
\begin{align*}
(u_{3},v_{13},v_{23},w_{3})  &  \in\lbrack\overline{B}(u_{2},\gamma
)\times\overline{B}(v_{12},\gamma)\times\overline{B}(v_{22},\gamma
)\times\overline{B}(w_{2},\gamma)]\cap C_{1},\\
(u_{4},v_{14},v_{24},w_{4})  &  \in\lbrack\overline{B}(u_{2},\gamma
)\times\overline{B}(v_{12},\gamma)\times\overline{B}(v_{22},\gamma
)\times\overline{B}(w_{2},\gamma)]\cap C_{2},\\
(u_{5},v_{15},v_{25},w_{5})  &  \in\lbrack\overline{B}(u_{2},\gamma
)\times\overline{B}(v_{12},\gamma)\times\overline{B}(v_{22},\gamma
)\times\overline{B}(w_{2},\gamma)]\cap C_{3}.
\end{align*}

In conclusion, there exist
\begin{align*}
w_{0}^{\ast}  &  \in S_{Z^{\ast}},\\
(u_{3}^{\ast},v_{13}^{\ast},0,0)  &  \in\widehat{N}(C_{1},(u_{3},v_{13}%
,v_{23},w_{3}))\iff u_{3}^{\ast}\in\widehat{D}^{\ast}F_{1}(u_{3}%
,v_{13})(-v_{13}^{\ast}),\\
(u_{4}^{\ast},0,v_{24}^{\ast},0)  &  \in\widehat{N}(C_{2},(u_{4},v_{14}%
,v_{24},w_{4}))\iff u_{4}^{\ast}\in\widehat{D}^{\ast}F_{2}(u_{4}%
,v_{14})(-v_{24}^{\ast}),\\
(u_{5}^{\ast},v_{15}^{\ast},v_{25}^{\ast},w_{5}^{\ast})  &  \in\overline
{B}_{X^{\ast}}\times\overline{B}_{Y_{1}^{\ast}}\times\overline{B}_{Y_{2}%
^{\ast}}\times\overline{B}_{Z^{\ast}}%
\end{align*}

\noindent such that%
\begin{align*}
&  (-u_{3}^{\ast}-u_{4}^{\ast}-(b+2\rho)u_{5}^{\ast},-v_{13}^{\ast}%
-\alpha(1+2\rho)v_{15}^{\ast},-v_{24}^{\ast}-\alpha(1+2\rho)v_{25}^{\ast
},-w_{0}^{\ast}-(b+2\rho)w_{5}^{\ast})\\
&  \in\widehat{N}(C_{3},(u_{5},v_{15},v_{25},w_{5})),
\end{align*}

\noindent i.e.%
\begin{align*}
-u_{3}^{\ast}-u_{4}^{\ast}-(b+2\rho)u_{5}^{\ast}  &  =0,\\
(-v_{13}^{\ast}-\alpha(1+2\rho)v_{15}^{\ast},-v_{24}^{\ast}-\alpha
(1+2\rho)v_{25}^{\ast})  &  \in\widehat{D}^{\ast}G(v_{15},v_{25},w_{5}%
)(w_{0}^{\ast}+(b+2\rho)w_{5}^{\ast}).
\end{align*}

Observe that%
\[
\left\Vert w_{0}^{\ast}+(b+2\rho)w_{5}^{\ast}\right\Vert \geq1-(b+2\rho
)>1-(b+2\varepsilon)>0,
\]

\noindent and denote%
\begin{align*}
x_{1}^{\ast}  &  :=\left\Vert w_{0}^{\ast}+(b+2\rho)w_{5}^{\ast}\right\Vert
^{-1}u_{3}^{\ast},\\
x_{2}^{\ast}  &  :=\left\Vert w_{0}^{\ast}+(b+2\rho)w_{5}^{\ast}\right\Vert
^{-1}u_{4}^{\ast},\\
t_{1}^{\ast}  &  :=-\left\Vert w_{0}^{\ast}+(b+2\rho)w_{5}^{\ast}\right\Vert
^{-1}v_{13}^{\ast},\\
t_{2}^{\ast}  &  :=-\left\Vert w_{0}^{\ast}+(b+2\rho)w_{5}^{\ast}\right\Vert
^{-1}v_{24}^{\ast},\\
w^{\ast}  &  :=\left\Vert w_{0}^{\ast}+(b+2\rho)w_{5}^{\ast}\right\Vert
^{-1}(w_{0}^{\ast}+(b+2\rho)w_{5}^{\ast}),\\
z_{1}^{\ast}  &  :=-\left\Vert w_{0}^{\ast}+(b+2\rho)w_{5}^{\ast}\right\Vert
^{-1}\alpha(1+2\rho)v_{15}^{\ast}\\
z_{2}^{\ast}  &  :=-\left\Vert w_{0}^{\ast}+(b+2\rho)w_{5}^{\ast}\right\Vert
^{-1}\alpha(1+2\rho)v_{25}^{\ast}.
\end{align*}

In conclusion, one has%
\begin{align*}
x_{1}^{\ast}  &  \in\widehat{D}^{\ast}F_{1}(u_{3},v_{13})(t_{1}^{\ast}),\\
x_{2}^{\ast}  &  \in\widehat{D}^{\ast}F_{2}(u_{4},v_{24})(t_{2}^{\ast}),\\
(t_{1}^{\ast}+z_{1}^{\ast},t_{2}^{\ast}+z_{2}^{\ast})  &  \in\widehat{D}%
^{\ast}G(v_{15},v_{25},w_{5})(w^{\ast}),
\end{align*}

\noindent where%
\begin{align}
\left\Vert w^{\ast}\right\Vert  &  =1,\label{e.1}\\
\left\Vert x_{1}^{\ast}+x_{2}^{\ast}\right\Vert  &  =\frac{\left\Vert
u_{3}^{\ast}+u_{4}^{\ast}\right\Vert }{\left\Vert w_{0}^{\ast}+(b+2\rho
)w_{5}^{\ast}\right\Vert }\leq\frac{b+2\rho}{1-(b+2\rho)}<c\label{e.2}\\
\left\Vert z_{1}^{\ast}\right\Vert  &  \leq\frac{\alpha(1+2\rho)}{1-(b+2\rho
)}<\delta,\label{e.3}\\
\left\Vert z_{2}^{\ast}\right\Vert  &  \leq\frac{\alpha(1+2\rho)}{1-(b+2\rho
)}<\delta. \label{e.4}%
\end{align}

Remark that $(u_{3},v_{13})\in\operatorname*{Gr}F_{1}$ and
\begin{align*}
(u_{3},v_{13})  &  \in\overline{B}(u_{2},\gamma)\times\overline{B}%
(y_{12},\gamma)\subset\overline{B}(x_{0},2\gamma)\times\overline{B}%
(y_{10},2\gamma)\\
&  \subset B(\widehat{x},2^{-1}r)\times B(\widehat{y}_{1},2^{-1}r)\subset
B(\overline{x},r)\times B(\overline{y}_{2},r).
\end{align*}

\noindent Analogously, $(u_{4},v_{24})\in\operatorname*{Gr}F_{2}\cap\lbrack
B(\overline{x},r)\times B(\overline{y}_{2},r)]$ and $(v_{15},v_{25},w_{5}%
)\in\operatorname*{Gr}G\cap\lbrack B(\overline{y}_{1},r)\times B(\overline
{y}_{2},r)\times B(\overline{z},r)].$ Using now (\ref{ineq_x12}) and
(\ref{e.2}), we get that%
\[
c\leq\left\Vert x_{1}^{\ast}+x_{2}^{\ast}\right\Vert <c,
\]

\noindent a contradiction. It follows that $v=z_{0}$ is the sole possible
situation and the conclusion follows.$\hfill\square$

\bigskip

We remark that the conclusion of Theorem \ref{coder_Fr} could be deduced (with
a different technique) using as an intermediate step some estimations for the
strong slope of the lower semicontinuous envelope $\varphi_{R},$ as done in
\cite{NTT}, \cite{DNS2011}. We prefered here this direct aproach for clarity
and for the sake of including our result in the framework opened by
Mordukhovich and Shao \cite{MordShao1996} and Penot \cite{Pen1998a}.

\bigskip

We want to discuss in the sequel the possible relations between the
coderivative conditions which assure the linear openness/metric regularity of
$H$ and the main result, i.e. Theorem \ref{mainresult}. To this, let us
formulate first an auxiliary result.

\begin{pr}
\label{estim_coder}(A) Let $X,Y$ be Banach spaces, $\Gamma:X\rightrightarrows
Y$ be a multifunction and $(\overline{x},\overline{y})\in\operatorname*{Gr}%
\Gamma.$

(i) If $\Gamma$ has the Aubin property around $(\overline{x},\overline{y}),$
then there exists $\alpha>0$ arbitrarily close to $\operatorname{lip}%
\Gamma(\overline{x},\overline{y})$ and $r>0$ such that, for every
$(x,y)\in\operatorname*{Gr}\Gamma\cap\lbrack B(\overline{x},r)\times
B(\overline{y},r)],$ $y^{\ast}\in Y^{\ast}$ and $x^{\ast}\in\widehat{D}^{\ast
}\Gamma(x,y)(y^{\ast}),$%
\[
\left\Vert x^{\ast}\right\Vert \leq\alpha\left\Vert y^{\ast}\right\Vert .
\]

(ii) If $\Gamma$ is open at linear rate around $(\overline{x},\overline{y}),$
then there exists $\alpha>0$ arbitrarily close to $\operatorname{lop}%
\Gamma(\overline{x},\overline{y})$ and $r>0$ such that, for every
$(x,y)\in\operatorname*{Gr}\Gamma\cap\lbrack B(\overline{x},r)\times
B(\overline{y},r)],$ $y^{\ast}\in Y^{\ast}$ and $x^{\ast}\in\widehat{D}^{\ast
}\Gamma(x,y)(y^{\ast}),$%
\[
\left\Vert x^{\ast}\right\Vert \geq\alpha\left\Vert y^{\ast}\right\Vert .
\]

(B) Let $X,Y,Z$ be Banach spaces, $\Phi:X\times Y\rightrightarrows Z$ be a
multifunction and $((\overline{x},\overline{y}),\overline{z})\in
\operatorname*{Gr}\Phi.$

(i) If $\Phi$ has the Aubin property with respect to $x,$ uniformly in $y$
around $((\overline{x},\overline{y}),\overline{z}),$ then there exists
$\alpha>0$ arbitrarily close to $\widehat{\operatorname{lip}}_{x}%
\Phi((\overline{x},\overline{y}),\overline{z})$ and $r>0$ such that, for every
$y\in B(\overline{y},r),$ every $(x,z)\in\operatorname*{Gr}\Phi_{y}\cap\lbrack
B(\overline{x},r)\times B(\overline{z},r)],$ $z^{\ast}\in Z^{\ast}$ and
$x^{\ast}\in\widehat{D}^{\ast}\Phi_{y}(x,z)(z^{\ast}),$%
\[
\left\Vert x^{\ast}\right\Vert \leq\alpha\left\Vert z^{\ast}\right\Vert .
\]

(ii) If $\Phi$ is open at linear rate with respect to $x,$ uniformly in $y$
around $((\overline{x},\overline{y}),\overline{z}),$ then there exists
$\alpha>0$ arbitrarily close to $\widehat{\operatorname{lop}}_{x}%
\Phi((\overline{x},\overline{y}),\overline{z})$ and $r>0$ such that, for every
$y\in B(\overline{y},r),$ every $(x,z)\in\operatorname*{Gr}\Phi_{y}\cap\lbrack
B(\overline{x},r)\times B(\overline{z},r)],$ $z^{\ast}\in Z^{\ast}$ and
$x^{\ast}\in\widehat{D}^{\ast}\Phi_{y}(x,z)(z^{\ast}),$%
\[
\left\Vert x^{\ast}\right\Vert \geq\alpha\left\Vert z^{\ast}\right\Vert .
\]

\end{pr}

\noindent\textbf{Proof.} The assertion from $(A)$ are well-known. See, for
instance, \cite[Theorems 1.43, 1.54]{Mor2006}. Concerning $(B),$ the proof is
very similar to the one of $(A)$ and it mimics the proofs of the mentioned
theorems. For the readers convenience, let us prove the item $(i)$ from $(B).$

Suppose that $\Phi$ has the Aubin property with respect to $x,$ uniformly in
$y$ around $((\overline{x},\overline{y}),\overline{z}),$ hence there exists
$\alpha>0$ arbitrarily close to $\widehat{\operatorname{lip}}_{x}%
\Phi((\overline{x},\overline{y}),\overline{z})$ and $\varepsilon>0$ such that,
for every $y\in B(\overline{y},\varepsilon)$ and every $x,u\in B(\overline
{x},\varepsilon),$%
\begin{equation}
\Phi_{y}(x)\cap B(\overline{z},\varepsilon)\subset\Phi_{y}(u)+\alpha\left\Vert
x-u\right\Vert B_{Z}. \label{Aub_phi}%
\end{equation}

Consider $r:=2^{-1}\varepsilon$ and take arbitrary $y\in B(\overline{y},r),$
$(x,z)\in\operatorname*{Gr}\Phi_{y}\cap\lbrack B(\overline{x},r)\times
B(\overline{z},r)],$ $z^{\ast}\in Z^{\ast},$ $x^{\ast}\in\widehat{D}^{\ast
}\Phi_{y}(x,z)(z^{\ast})$ and $\delta>0.$ Employing the definition of the
Fr\'{e}chet coderivatives, one may find $\gamma<\min\{2^{-1}\varepsilon
,2^{-1}\alpha\varepsilon\}$ such that%
\begin{equation}
\left\langle x^{\ast},u-x\right\rangle -\left\langle z^{\ast},w-z\right\rangle
\leq\delta(\left\Vert u-x\right\Vert +\left\Vert w-z\right\Vert )
\label{coder_phi}%
\end{equation}

\noindent for every $(u,w)\in\operatorname*{Gr}\Phi_{y}\cap\lbrack
B(x,\alpha^{-1}\gamma)\times B(z,\gamma)].$

Take now arbitrary $u\in B(x,\alpha^{-1}\gamma).$ Then%
\[
\left\Vert u-\overline{x}\right\Vert \leq\left\Vert u-x\right\Vert +\left\Vert
x-\overline{x}\right\Vert \leq\alpha^{-1}\gamma+r<\varepsilon,
\]

\noindent hence one may apply (\ref{Aub_phi}) for this $u$ and $z\in\Phi
_{y}(x)\cap B(\overline{z},\varepsilon)$ in order to find $w\in\Phi_{y}(u)$
such that%
\[
\left\Vert w-u\right\Vert \leq\alpha\left\Vert x-u\right\Vert \leq\gamma.
\]

Consequently, using (\ref{coder_phi}), one obtains%
\[
\alpha^{-1}\gamma\left\Vert x^{\ast}\right\Vert \leq\gamma\left\Vert z^{\ast
}\right\Vert +\delta\gamma(1+\alpha^{-1}).
\]

\noindent As $\delta$ is arbitrary positive, the conclusion follows.

The proof of the item $(ii)$ from $(B)$ takes into account two facts. The
first one is the obvious equivalence%
\[
x^{\ast}\in\widehat{D}^{\ast}\Phi_{y}(x,z)(z^{\ast})\iff-z^{\ast}\in
\widehat{D}^{\ast}\Phi_{y}^{-1}(z,x)(-x^{\ast}).
\]

\noindent The second one is more technical, and it concerns a link between the
regularity notions in the parametric case. More precisely, following the lines
of the proof of Theorem \ref{link_around}, if $\Phi$ is open at linear rate
with respect to $x,$ uniformly in $y$ around $((\overline{x},\overline
{y}),\overline{z}),$ then $\Phi_{y}^{-1}$ is Aubin around $(\overline
{z},\overline{x}),$ for every $y$ in a neighborhood which can be taken exactly
the same as in the definition of the uniform openness in $y.$ Then $(A),$ item
$(i),$ can be applied and the conclusion follows.$\hfill\square$

\bigskip

Let us make the final step of our plan to derive again Theorem
\ref{mainresult} on Asplund spaces, by means of coderivative conditions in
Theorem \ref{coder_Fr}.

\begin{thm}
Let $X,Y_{1},Y_{2}$ and $Z$ be Asplund spaces, $F_{1}:X\rightrightarrows
Y_{1},$ $F_{2}:X\rightrightarrows Y_{2},$ $G:Y_{1}\times Y_{2}%
\rightrightarrows Z$ be closed-graph multifunctions, and $(\overline
{x},\overline{y}_{1},\overline{y}_{2},\overline{z})\in X\times Y_{1}\times
Y_{2}\times Z$ be such that $\overline{z}\in G(\overline{y}_{1},\overline
{y}_{2}),(\overline{y}_{1},\overline{y}_{2})\in F_{1}(\overline{x})\times
F_{2}(\overline{x}).$ Suppose that the assumptions (i)-(v) of Theorem
\ref{mainresult} are satisfied. Then the sets $C_{1},C_{2},C_{3}$ given by
(\ref{C13}) are allied at $(\overline{x},\overline{y}_{1},\overline{y}%
_{2},\overline{z}).$ Moreover, there exist $L,M,C,D>0$ arbitrarily close to
$\operatorname*{reg}F_{1}(\overline{x},\overline{y}_{1}),$
$\operatorname*{lip}F_{2}(\overline{x},\overline{y}_{2}),$ $\widehat
{\operatorname*{reg}}_{y_{1}}G((\overline{y}_{1},\overline{y}_{2}%
),\overline{z}),$ $\widehat{\operatorname*{lip}}_{y_{2}}G((\overline{y}%
_{1},\overline{y}_{2}),\overline{z})$ such that $LMCD<1$ and for every
$\varepsilon>0$ arbitrary small, the relation \ref{C} is satisfied with
\[
c:=\frac{1}{LC}-MD-\varepsilon>0.
\]

In conclusion, for every $a\in(0,c),$ $H$ is metrically regular at
$(\overline{x},\overline{z})$ with constant $a.$ If, moreover, $(F_{1}%
,F_{2}),G$ are locally composition-stable around $(\overline{x},(\overline
{y}_{1},\overline{y}_{2}),\overline{z}),$ then for every $a\in(0,c),$ $H$ is
metrically regular around $(\overline{x},\overline{z})$ with constant $a.$ As
consequence,%
\[
\operatorname*{reg}H(\overline{x},\overline{z})\leq\rho_{0},
\]

\noindent where $\rho_{0}$ is given by (\ref{ro0}).
\end{thm}

\noindent\textbf{Proof.} Suppose that the multifunctions $F_{1},F_{2},G$
satisfy the assumptions (ii)-(iv) of Theorem \ref{mainresult} and consider
arbitrary sequences $(x_{n},y_{1n})\overset{\operatorname*{Gr}F_{1}%
}{\longrightarrow}(\overline{x},\overline{y}_{1}),$ $(u_{n},v_{2n}%
)\overset{\operatorname*{Gr}F_{2}}{\longrightarrow}(\overline{x},\overline
{y}_{2}),$ $(b_{1n},b_{2n},c_{n})\overset{\operatorname*{Gr}G}{\longrightarrow
}(\overline{y}_{1},\overline{y}_{2},\overline{z})$ and $x_{n}^{\ast}%
\in\widehat{D}^{\ast}F_{1}(x_{n},y_{1n})(y_{1n}^{\ast}),$ $u_{n}^{\ast}%
\in\widehat{D}^{\ast}F_{2}(u_{n},v_{2n})(v_{2n}^{\ast}),$ $(b_{1n}^{\ast
},b_{2n}^{\ast})\in\widehat{D}^{\ast}G(b_{1n},b_{2n},c_{n})(c_{n}^{\ast})$
such that the relations $(x_{n}^{\ast}+u_{n}^{\ast})\rightarrow0,$ $\left(
y_{1n}^{\ast}+b_{1n}^{\ast}\right)  \rightarrow0,\left(  v_{2n}^{\ast}%
+b_{2n}^{\ast}\right)  \rightarrow0,\left(  c_{n}^{\ast}\right)  \rightarrow0$
hold. Then for every $n$, $b_{1n}^{\ast}\in\widehat{D}^{\ast}G_{b_{2n}}%
(b_{1n},c_{n})(c_{n}^{\ast})$ and $b_{2n}^{\ast}\in\widehat{D}^{\ast}%
G_{b_{1n}}(b_{2n},c_{n})(c_{n}^{\ast}).$ Using now Theorem \ref{estim_coder},
there exist $L,M,C,D>0$ arbitrarily close to $\operatorname*{reg}%
F_{1}(\overline{x},\overline{y}_{1}),$ $\operatorname*{lip}F_{2}(\overline
{x},\overline{y}_{2}),$ $\widehat{\operatorname*{reg}}_{y_{1}}G((\overline
{y}_{1},\overline{y}_{2}),\overline{z}),$ $\widehat{\operatorname*{lip}%
}_{y_{2}}G((\overline{y}_{1},\overline{y}_{2}),\overline{z})$ such that
$LMCD<1$ and, for every $n,$ one has%
\begin{align}
\left\Vert x_{n}^{\ast}\right\Vert  &  \geq L^{-1}\left\Vert y_{1n}^{\ast
}\right\Vert ,\nonumber\\
\left\Vert u_{n}^{\ast}\right\Vert  &  \leq M\left\Vert v_{2n}^{\ast
}\right\Vert ,\label{allied_est}\\
\left\Vert b_{1n}^{\ast}\right\Vert  &  \geq C^{-1}\left\Vert c_{n}^{\ast
}\right\Vert ,\nonumber\\
\left\Vert b_{2n}^{\ast}\right\Vert  &  \leq D\left\Vert c_{n}^{\ast
}\right\Vert .\nonumber
\end{align}

As $\left(  c_{n}^{\ast}\right)  \rightarrow0,$ it follows successively from
(\ref{allied_est}) and the assumptions made that $\left(  b_{2n}^{\ast
}\right)  \rightarrow0,$ $\left(  v_{2n}^{\ast}\right)  \rightarrow0,$
$(u_{n}^{\ast})\rightarrow0,$ $(x_{n}^{\ast})\rightarrow0,$ $(y_{1n}^{\ast
})\rightarrow0,$ $\left(  b_{1n}^{\ast}\right)  \rightarrow0.$ In conclusion,
the sets $C_{1},C_{2},C_{3}$ are allied at $(\overline{x},\overline{y}%
_{1},\overline{y}_{2},\overline{z}).$

Using again the properties (i)-(iv) and Theorem \ref{estim_coder}, one gets
the existence of $r>0$ and of $L,M,C,D>0$ as above such that:

\begin{enumerate}
\item for every $(x,y_{1})\in\operatorname*{Gr}F_{1}\cap\lbrack B(\overline
{x},r)\times B(\overline{y}_{1},r)],$ $y_{1}^{\ast}\in Y^{\ast}$ and $x^{\ast
}\in\widehat{D}^{\ast}F_{1}(x,y_{1})(y_{1}^{\ast}),$%
\begin{equation}
\left\Vert x^{\ast}\right\Vert \geq L^{-1}\left\Vert y_{1}^{\ast}\right\Vert ;
\label{r1}%
\end{equation}

\item for every $(x,y_{2})\in\operatorname*{Gr}F_{2}\cap\lbrack B(\overline
{x},r)\times B(\overline{y}_{2},r)],$ $y_{2}^{\ast}\in Y^{\ast}$ and $x^{\ast
}\in\widehat{D}^{\ast}F_{2}(x,y_{2})(y_{2}^{\ast}),$%
\begin{equation}
\left\Vert x^{\ast}\right\Vert \leq M\left\Vert y_{2}^{\ast}\right\Vert ;
\label{r2}%
\end{equation}

\item for every $y_{2}\in B(\overline{y}_{2},r),$ every $(y_{1},z)\in
\operatorname*{Gr}G_{y_{2}}\cap\lbrack B(\overline{y}_{1},r)\times
B(\overline{z},r)],$ $z^{\ast}\in Z^{\ast}$ and $y_{1}^{\ast}\in\widehat
{D}^{\ast}G_{y_{2}}(y_{1},z)(z^{\ast}),$%
\begin{equation}
\left\Vert y_{1}^{\ast}\right\Vert \geq C^{-1}\left\Vert z^{\ast}\right\Vert ;
\label{r3}%
\end{equation}

\item for every $y_{1}\in B(\overline{y}_{1},r),$ every $(y_{2},z)\in
\operatorname*{Gr}G_{y_{1}}\cap\lbrack B(\overline{y}_{2},r)\times
B(\overline{z},r)],$ $z^{\ast}\in Z^{\ast}$ and $y_{2}^{\ast}\in\widehat
{D}^{\ast}G_{y_{1}}(y_{2},z)(z^{\ast}),$%
\begin{equation}
\left\Vert y_{2}^{\ast}\right\Vert \leq D\left\Vert z^{\ast}\right\Vert .
\label{r4}%
\end{equation}

\end{enumerate}

Choose next arbitrary $\varepsilon>0$ and $\delta\in(0,r)$ such that
$(L^{-1}+M)\delta<\varepsilon.$ Take ${(u_{1},v_{1})\in}\operatorname*{Gr}%
F_{1}\cap\lbrack B(\overline{x},r)\times B(\overline{y}_{1},r)],$
${(u_{2},v_{2})\in}\operatorname*{Gr}F_{2}\cap\lbrack B(\overline{x},r)\times
B(\overline{y}_{2},r)],$ $(t_{1},t_{2},w)\in\operatorname*{Gr}G\cap\lbrack
B(\overline{y}_{1},r)\times B(\overline{y}_{2},r)\times B(\overline{z},r)],$
$w^{\ast}\in S_{Z^{\ast}},$ $t_{1}^{\ast}\in Y_{1}^{\ast},$ $t_{2}^{\ast}\in
Y_{2}^{\ast},$ $z_{1}^{\ast}\in\delta B_{Y_{1}^{\ast}},$ $z_{2}^{\ast}%
\in\delta B_{Y_{2}^{\ast}}$ such that $(z_{1}^{\ast}+t_{1}^{\ast},z_{2}^{\ast
}+t_{2}^{\ast})\in\widehat{D}^{\ast}G((t_{1},t_{2}),w)(w^{\ast})$ and
$x_{1}^{\ast}\in\widehat{D}^{\ast}F_{1}(u_{1},v_{1})(t_{1}^{\ast}),$
$x_{2}^{\ast}\in\widehat{D}^{\ast}F_{2}(u_{2},v_{2})(t_{2}^{\ast}).$ Then
$z_{1}^{\ast}+t_{1}^{\ast}\in\widehat{D}^{\ast}G_{t_{2}}(t_{1},w)(w^{\ast})$
and $z_{2}^{\ast}+t_{2}^{\ast}\in\widehat{D}^{\ast}G_{t_{1}}(t_{2},w)(w^{\ast
}),$ and using (\ref{r1})-(\ref{r4}), one gets%
\begin{align*}
\left\Vert x_{1}^{\ast}+x_{2}^{\ast}\right\Vert  &  \geq\left\Vert x_{1}%
^{\ast}\right\Vert -\left\Vert x_{2}^{\ast}\right\Vert \geq L^{-1}\left\Vert
t_{1}^{\ast}\right\Vert -M\left\Vert t_{2}^{\ast}\right\Vert \\
&  \geq L^{-1}\left(  \left\Vert z_{1}^{\ast}+t_{1}^{\ast}\right\Vert
-\left\Vert z_{1}^{\ast}\right\Vert \right)  -M(\left\Vert z_{2}^{\ast}%
+t_{2}^{\ast}\right\Vert +\left\Vert z_{2}^{\ast}\right\Vert )\\
&  \geq L^{-1}C^{-1}\left\Vert w^{\ast}\right\Vert -MD\left\Vert w^{\ast
}\right\Vert -(L^{-1}+M)\delta\\
&  >L^{-1}C^{-1}-MD-\varepsilon.
\end{align*}

It follows that relation \ref{C} is satisfied for $c=L^{-1}C^{-1}%
-MD-\varepsilon>0,$ i.e. the conclusion follows.$\hfill\square$


\bigskip

\begin{thebibliography}{99}                                                                                               %


\bibitem {Arut2007}A.V. Arutyunov, \textit{Covering mapping in metric spaces,
and fixed points}, Doklady Mathematics, 76 (2007), 665--668.

\bibitem {Arut2009}A.V. Arutyunov, \textit{Stability of coincidence points and
properties of covering mappings,} Mathematical Notes, 86 (2009), 153--158.

\bibitem {AAGDO}A. Arutyunov, E. Avakov, B. Gel'man, A. Dmitruk, V.
Obukhovskii, \textit{Locally covering maps in metric spaces and coincidence
points,} Journal of Fixed Points Theory and Applications, 5 (2009), 105--127.

\bibitem {RefAz}D. Az\'{e}, \textit{A unified theory for metric regularity of
multifunctions}, Journal of Convex Analysis, 13 (2006), 225--252.

\bibitem {RefBonS}J.F. Bonnans, A. Shapiro, \textit{Perturbation Analysis of
Optimization Problems,} Springer-Verlag, New York, 2000.

\bibitem {BD}J.M. Borwein, A.L. Dontchev, \textit{On the Bartle-Graves
theorem}, Proceedings of the American Mathematical Society, 131 (2003), 2553--2560.

\bibitem {[10]}J.M. Borwein, Q.J. Zhu, \textit{Techniques of Variational
Analysis,} CMS Books in Mathematics, 20, Springer-Verlag, New York, 2005.

\bibitem {RefBorZu}J.M. Borwein, Q.J. Zhu, \textit{Viscosity solutions and
viscosity subderivatives in smooth Banach spaces with applications to metric
regularity}, SIAM Journal on Control and Optimization,\textbf{ }34 (1996), 1568--1591.

\bibitem {RefBorZ}J.M. Borwein, D.M. Zhuang, \textit{Verifiable necessary and
sufficient conditions for openness and regularity of set-valued maps}, Journal
of Mathematical Analysis and Applications, 134 (1988), 441--459.

\bibitem {[17]}F.H. Clarke, \textit{Optimization and Nonsmooth Analysis},
Canadian Mathematical Society Series of Monographs and Advanced Texts. John
Wiley \& Sons Inc., New York, 1983.

\bibitem {RefCom}R. Cominetti, \textit{Metric regularity, tangent cones, and
second-order optimality conditions}, Applied Mathematics \& Optimization, 21
(1990), 265--287.

\bibitem {DMO}A.V. Dmitruk, A.A. Milyutin, N.P. Osmolovskii,
\textit{Lyusternik's theorem and the theory of extrema,} Uspekhi Mat. Nauk, 35
(1980), 11--51.

\bibitem {dontchev}A.L. Dontchev, \textit{The Graves theorem revisited,}
Journal of Convex Analysis, 3 (1996), 45--53.

\bibitem {DonFra2010}A.L. Dontchev, H. Frankowska, \textit{Lyusternik-Graves
theorem and fixed points}, Proceedings of the American Mathematical Society,
139 (2011), 521--534.

\bibitem {DonFra2}A.L. Dontchev, H. Frankowska, \textit{Lyusternik-Graves
theorem and fixed points II}, Journal of Convex Analysis, appeared online.

\bibitem {DL}A.L. Dontchev, A.S. Lewis, \textit{Perturbations and metric
regularity}, Set-Valued Analysis, 13 (2005), 417--438.

\bibitem {DLR}A.L. Dontchev, A.S. Lewis, R.T. Rockafellar, \textit{The radius
of metric regularity}, Transactions of the American Mathematical Society, 355
(2003), 493--517.

\bibitem {Don-Quin-Zla}A.L. Dontchev, M. Quincampoix, N. Zlateva,
\textit{Aubin criterion for metric regularity}, Journal of Convex Analysis, 13
(2006), 281--297.

\bibitem {DontRock2009b}A.L. Dontchev, R.T. Rockafellar, \textit{Implicit
functions and solution mappings}, Springer, Berlin, 2009.

\bibitem {DNS2011}M. Durea, H.T. Nguyen, R. Strugariu, \textit{Metric
regularity of epigraphical multivalued mappings and applications to vector
optimization}, Mathematical Programming, Serie B, accepted.

\bibitem {DurStr1}M. Durea, R. Strugariu, \textit{On some Fermat rules for
set-valued optimization problems}, Optimization, 60 (2011), 575--591.

\bibitem {DurStr4}M. Durea, R. Strugariu, \textit{Openness stability and
implicit multifunction theorems: Applications to variational systems},
Nonlinear Analysis, Theory, Methods and Applications, 75 (2012), 1246--1259.

\bibitem {DurStr5}M. Durea, R. Strugariu, \textit{Chain rules for linear
openness in general Banach spaces}, SIAM Journal on Optimization, 22 (2012), 899--913.

\bibitem {DurStr2012}M. Durea, R. Strugariu, \textit{Chain rules for linear
openness in metric spaces and applications}, Mathematical Programming, Serie
A, DOI: 10.1007/s10107-012-0598-8.

\bibitem {Gra}L.M. Graves,\textit{ Some mapping theorems}, Duke Math. Journal,
17 (1950), 111--114.

\bibitem {RefIo3}A. Ioffe, \textit{Regular points of Lipschitz functions},
Transactions of the American Mathematical Society, 251 (1979), 61--69.

\bibitem {RefIo}A.D. Ioffe, \textit{On perturbation stability of metric
regularity}, Set-Valued Analysis, 9 (2001), 101--109.

\bibitem {I Cont}A.D. Ioffe, \textit{On regularity estimates for mappings
between embedded manifolds}, Control and Cybernetics, 36 (2007), 659--668.

\bibitem {RefIoffe}A.D. Ioffe, \textit{On regularity concepts in variational
analysis}, Fixed Point Theory and Applications, 8 (2010), 339--363.

\bibitem {Ioffe2000}A.D. Ioffe, \textit{Metric regularity and subdifferential
calculus}, Uspekhi Mat. Nauk, 55 (2000), no. 3 (333), 103--162; English
translation in Math. Surveys, 55 (2000), 501--558.

\bibitem {Ioffe2010b}A.D. Ioffe, \textit{Towards variational analysis in
metric spaces: metric regularity and fixed points}, Mathematical Programming,
Serie B, 123 (2010), 241--252.

\bibitem {SiamReFi}A.D. Ioffe, \textit{Regularity on a Fixed Set}, SIAM
Journal on Optimization, \textbf{21}, No. 4, (2011), 1345-1370.

\bibitem {RefJAus}A. Jourani, \textit{On metric regularity of multifunctions},
Bulletin of the Australian Mathematical Society, 44 (1991).

\bibitem {Jo-Th1}A. Jourani, L. Thibault, \textit{Approximate subdifferential
and metric regularity : Finite dimensional case}, Mathematical Programming, 47
(1990), 203--218.

\bibitem {Jo-Th}A. Jourani, L. Thibault, \textit{Approximations and metric
regularity in mathematical programming in Banach spaces}, Mathematics of
Operations Research, 18 (1993), 390--401.

\bibitem {JouT1}A. Jourani, L. Thibault, \textit{Metric regularity for
strongly compactly Lipschitzian mappings}, Nonlinear Analysis TMA, 24 (1995), 229--240.

\bibitem {1JT2}A. Jourani, L. Thibault, \textit{Verifiable conditions for
openness and metric regularity of multivalued mappings in Banach spaces},
Transactions of the American Mathematical Society, 347 (1995), 1255--1268.

\bibitem {RefJT}A. Jourani, L. Thibault, \textit{Metric inequality and
subdifferential calculus in Banach spaces}, Set-Valued Analysis, 3 (1995), 87--100.

\bibitem {JT1}A. Jourani, L. Thibault, \textit{Coderivatives of multivalued
mappings, locally compact cones and metric regularity}, Nonlinear Analysis,
TMA, 35 (1994), 925--945.

\bibitem {AK}A.Y. Kruger, \textit{A covering theorem for set-valued mappings,}
Optimization, 19 (1988), 763--780.

\bibitem {[18]}A.V. Dmitruk, A.Y. Kruger, \textit{Metric regularity and
systems of generalized equations}, Journal of Mathematical Analysis and
Applications, 342 (2008), 864--873.

\bibitem {LPX2011}S. Li, J.-P. Penot, X. Xue, \textit{Codifferential calculus,
}Set-Valued and Variational Analysis, 19 (2011), 505--536.

\bibitem {Lyu}L.A. Lyusternik, \textit{On the conditional extrema of
functionals}, Mat. Sbornik 41 (1934) 390--401.

\bibitem {Mor1980}B. S. Mordukhovich, \textit{Metric approximations and
necessary optimality conditions for general classes of extremal problems,
}Soviet. Math. Dokl., 22 (1980), 526--530.

\bibitem {Mor2006}B.S. Mordukhovich, \textit{Variational Analysis and
Generalized Differentiation}, Vol. I: Basic Theory, Vol. II: Applications,
Springer, Grundlehren der mathematischen Wissenschaften (A Series of
Comprehensive Studies in Mathematics), Vol. 330 and 331, Berlin, 2006.

\bibitem {MordShao1996}B. S. Mordukhovich, Y. Shao, \textit{Nonconvex
differential calculus for infinite dimensional multifunctions}, Set-Valued
Analysis, 4 (1996), 205--336.

\bibitem {RefMorS}B.S. Mordukhovich, Y. Shao, \textit{Stability of set-valued
mappings in infinite dimensions: point criteria and applications}, SIAM
Journal on Control and Optimization, 35 (1997), 285--314.

\bibitem {RefMorW}B.S. Mordukhovich, B. Wang, \textit{Restrictive metric
regularity and generalized differential calculus in Banach spaces},
International Journal of Mathematical Sciences, 50 (2004), 2650--2683.

\bibitem {NT2001}H.V. Ngai, M. Th\'{e}ra, \textit{Metric inequality,
subdifferential calculus and applications,} Set-Valued Analysis, 9 (2001), 187--216.

\bibitem {RefNT3}H.V. Ngai, M. Th\'{e}ra, \textit{Error bounds and implicit
multifunctions in smooth Banach spaces and applications to optimization},
Set-Valued Analysis,\emph{ }12 (2004), 195--223.

\bibitem {NT2008}H.V. Ngai, M. Th\'{e}ra, \textit{Error bounds in metric
spaces and application to the perturbation stability of metric regularity},
SIAM Journal on Optimization, 19 (2008), 1--20.

\bibitem {NTT}H.V. Ngai, H.T. Nguyen, M. Th\'{e}ra, \textit{Implicit
multifunction theorems in complete metric spaces, } Mathematical Programming,
Serie B, accepted.

\bibitem {NTT2}H.V. Ngai, H.T. Nguyen, M. Th\'{e}ra, \textit{Metric regularity
of the sum of multifunctions and applications}, available at http://www.optimization-online.org/DB\_HTML/2011/12/3291.html.x

\bibitem {[69]}J. Outrata, M. Kocvara, J. Zowe, \textit{Nonsmooth Approach to
Optimization Problems with Equilibrium Constraints}, Nonconvex Optimization
and Its Applications, 28, Kluwer Academic Publishers, Dordrecht, 1998.

\bibitem {Pen1998}J.-P. Penot, \textit{Cooperative behavior of functions, sets
and relations,} Mathematical Methods of Operations Research, 48 (1998), 229--246.

\bibitem {Pen1998a}J.-P. Penot, \textit{Compactness properties, openness
criteria and coderivatives}, Set-Valued Analysis, 6 (1998), 363--380.

\bibitem {RefPe}J.-P. Penot, \textit{Regularity, openess and Lipschitzian
behavior of multifunctions}, Nonlinear Analysis, 13 (1989), 629--643.

\bibitem {Rob1976}S.M. Robinson, \textit{Regularity and stability for convex
multivalued functions}, Mathematics of Operations Research, 1 (1976), 130--143.

\bibitem {[74]}S.M. Robinson, \textit{Strongly regular generalized equations},
Mathematics of Operations Research, 5 (1980), 43--62.

\bibitem {RocWet}R.T. Rockafellar, R. Wets, \textit{Variational Analysis},
Springer, Grundlehren der mathematischen Wissenschaften (A Series of
Comprehensive Studies in Mathematics), Vol. 317, Berlin, 1998.

\bibitem {Schiro}W. Schirotzek, \textit{Nonsmooth Analysis}, Springer, Berlin, (1998).

\bibitem {Urs1975}C. Ursescu, \textit{Multifunctions with closed convex
graph}, Czech. Math. J. 25 (1975), 438--441.

\bibitem {Urs1996}C. Ursescu, \textit{Inherited openness}, Revue Roumaine des
Math\'{e}matiques Pures et Appliqu\'{e}es, 41 (1996), 5--6, 401--416.

\bibitem {Zng}X.Y. Zheng, K.F. Ng, \textit{Metric regularity and constraint
qualifications for convex inequalities on Banach spaces}, SIAM Journal on
Optimization, 14 (2003), 757--772.
\end{thebibliography}
\end{document}